\documentclass[final,onethmnum, onefignum, onetabnum]{siamltex}
\usepackage{amsmath}
\usepackage{amssymb}
\usepackage[lined,boxruled,linesnumbered]{algorithm2e}
\usepackage{graphicx}
\usepackage{subfig}

\usepackage{pgf}
\usepackage{tikz}
\usepackage{color}

\newcommand{\tenvec}{\operatorname{vec}}
\usetikzlibrary{positioning,arrows,automata,fit,shapes}

\title{Multigrid methods for tensor structured Markov chains with low rank approximation} 

\author{Matthias Bolten$^\ast$, Karsten Kahl$^\ast$ and Sonja Sokolovi\'{c}$^\ast$}

\begin{document}
	\maketitle
	
	\renewcommand{\thefootnote}{\fnsymbol{footnote}}
	\newcommand{\FF}{F\hspace{-0.075cm}F}
	\footnotetext[1]{Fachbereich C, Mathematik und Naturwissenschaften, Bergische Universit\"at Wuppertal, 42097 Wuppertal, Germany, \texttt{\{bolten,kkahl,sokolovic\}@math.uni-wuppertal.de}}
	\begin{abstract}
		Tensor structured Markov chains are part of stochastic models of many practical applications, e.g., 
		in the description of complex production or telephone networks. The most interesting question in 
		Markov chain models is the determination of the stationary distribution as a description of the 
		long term behavior of the system. This involves the computation of the eigenvector corresponding 
		to the dominant eigenvalue or equivalently the solution of a singular linear system of equations. 
		Due to the tensor structure of the models the dimension of the operators grows rapidly and a direct 
		solution without exploiting the tensor structure becomes infeasible. Algebraic multigrid methods have 
		proven to be efficient when dealing with Markov chains without using tensor structure. In this work we 
		present an approach to adapt the algebraic multigrid framework to the tensor frame, not only using the tensor 
		structure in matrix-vector multiplications, but also tensor structured coarse-grid operators and tensor representations 
		of the solution vector.
	\end{abstract}
	
	\begin{keywords}
		multigrid method, tensor truncation, Tensor Train, Markov chains, singular linear system
	\end{keywords}
	
	\begin{AMS}
		65F10, 65F50, 60J22, 65N55
	\end{AMS}
	
	\pagestyle{myheadings}
	\thispagestyle{plain}
	\markboth{M. BOLTEN, K. KAHL AND S. SOKOLOVI\'{C}}{MULTIGRID FOR TENSOR-STRUCTURED MARKOV CHAINS}
	
	\section{Introduction}
	Consider the transition rate matrix $A$ of an irreducible con\-tin\-u\-ous--time Markov chain and the task to find a state vector $x$ such that
	\begin{equation}\label{eq:Ax0}
		Ax=0.
	\end{equation}
	where $x \ge 0$ component wise and the sum of the components is equal to one, i.e., $x$ can be interpreted as a vector of probabilities.
	
	The matrix $A \in \mathbb{R}^{n \times n}$ is singular with corresponding rank $n-1$ and has column sum zero, i.e., $\mathbf{1}^TA=0$ where $\mathbf{1}=[1, \cdots, 1]^T$. The solution of \eqref{eq:Ax0} is called the steady state vector and its existence and uniqueness (up to a scalar factor) is proven in~\cite{BermanPlemmons1994}.
	Many continuous--time Markov chains arising in practical applications are known for the so called state space explosion, see, e.g.,~\cite{BuchholzDayar2007}, especially those where the generator matrix $A$ has tensor structure 
	\begin{equation}\label{eq:A}
		A = \sum\limits_{t=1}^T\bigotimes_{j = 1}^J E_j^t, 
	\end{equation}
	where the matrices $E_j^j \in \mathbb{R}^{n_{j}\times n_{j}}$ describe 
	local transitions in each of the $J$ submodels and the other matrices $E_j^t \in \mathbb{R}^{n_{j}\times n_{j}}$ the synchronized transitions between the submodels. The sum over the Kronecker products of the matrices gives us the generator matrix of the whole system; cf.~\cite{Buchholz2000,BuchholzDayar2007}. Models of this type appear, e.g., in queuing theory~\cite{Chan1987,Chan1988,Kaufman1983} and analysis of stochastic automata networks~\cite{LangvilleStewart2004, PlateauStewart1997}.
	Note that most of the $E_j^t, t \neq j$ are the identity matrix and that an increase of the number of synchronized transitions results in an increase of the number of summands in \eqref{eq:Ax0}.
	In addition, the dimension of $A$ grows exponentially in the number of submodels $J$. Even for moderately large $J$ the storage and computational complexity of forming $A$ explicitly is prohibitive. 
	Thus, in order to be able to do any computations for such models, all operations with $A$ have to make use 
	of the compact format \eqref{eq:A}. Due to the fact that $x\in \mathbb{R}^{\prod_{j=1}^{J} n_{j}}$ it is again difficult to calculate the exact solution due to memory constraints.
	If the solution $x$ of \eqref{eq:Ax0} can be approximated by a vector $\tilde{x}$ of the form 
	\begin{equation}\label{tildex}
		\tilde{x} = \sum\limits_{i = 1}^R \bigotimes\limits_{j=1}^J x^i_j
	\end{equation} 
	with small $R$, the memory requirement for $x$ is reduced from $\prod_{j=1}^{J} n_{j}$ to $R \sum_{j=1}^{J}n_{j}$. In addition, using this format matrix-vector multiplications can now be calculated with even smaller computational cost.
	While this assumption is certainly fulfilled for the extreme case of non-interacting networks, see Theorem~\ref{the:noninteracting}, it is unclear whether this can be expected to be true in all applications of interest. Numerical experiments shown in section~\ref{sec:numtest} suggest that this is the case at least for the Markov chains considered here, c.f.\ Fig.~\ref{fig:rankplot}.
	
	The structure of $\tilde{x}$ can be interpreted as an outer vector product, which is equivalent to a $J$-way tensor in the canonical format with tensor rank $R$; cf.~\cite{Kiers2000, KoldaBader2009}. The approximation quality of $\tilde{x}$ depends on its tensor representation; cf.~\cite{Hackbusch2012}. 
	In this work we focus on the Tensor Train (TT) format~\cite{Oseledets2010,Oseledets2011,OseledetsDolgov2012} and its performance in a multigrid setting. Recent publications which deal with multigrid methods for Markov chains \cite{BoltenBrandtBrannickFrommerKahlLivshits2011,BrannickKahlSokolovic2014,BuchholzDayar2004,DeSterckManteuffelMcCormickMillerRugeSanders2010,DeSterckManteuffelMcCormickNguyenRuge2008,DeSterckManteuffelMcCormickRugeMillerPearsonSanders2010} do not discuss the use of tensors formats in a multigrid context, but show that for simple structured Markov chains a multigrid ansatz often works efficiently. Recently a multigrid approach has been proposed in \cite{DeSterckMiller2013} to obtain a low rank approximation to speed up the ALS algorithm.
	
	In the last years there have been two main directions in the development of multigrid or multilevel algorithms for Markov chains. On the one hand there are methods based on smoothed aggregation multigrid~\cite{DeSterckManteuffelMcCormickMillerRugeSanders2010,DeSterckManteuffelMcCormickNguyenRuge2008,DeSterckManteuffelMcCormickRugeMillerPearsonSanders2010}, on the other hand the bootstrap algebraic multigrid framework has also been investigated in the context of computing the stationary distribution~\cite{BoltenBrandtBrannickFrommerKahlLivshits2011,BrannickKahlSokolovic2014}. In addition in~\cite{BuchholzDayar2004} 
	multigrid methods for structured Markov chains like the ones considered in this paper are studied. They are based on aggregation of the submodels and do not maintain the structure of \eqref{eq:A}. Apart from this multigrid approach there are publications which apply to problems with generator matrices of the form \eqref{eq:A} which have a so called product form solution that can be then computed efficiently~\cite{Buchholz2008}. Note that one can interpret our approximation \eqref{tildex} of \eqref{eq:Ax0} as a product form solution if the number of summands is equal to one. Therefore it is obvious to interpret our approach as a combination and generalization of these two established solution techniques. Another recent approach for approximating the stationary distribution of tensor structured Markov chains by a low rank tensor was presented in~\cite{KressnerMacedo2014}, based on similar techniques for eigenvalue computations from~\cite{KressnerSteinlechnerUschmajew2013}. This approach does, however, not use multigrid techniques. 
	
	Tensor techniques have been used in multigrid methods in different settings.
	In \cite{BoermHiptmair2001} a tensorized representation of multigrid methods similar to the one presented here was used to construct robust multigrid methods for partial differential equations, including singular perturbations.
	Later a multigrid method to solve large scale Sylvester equations including a low rank representation of the current approximate solution was presented in \cite{GrasedyckHackbusch2007}.
	A related approach for the Lyapunov equation has been analyzed in \cite{VandereyckenVandewalle2009}.
	
	The remainder of this work is structured as follows. In sections \ref{sec:multigridbasics} and \ref{sec:tensorbasics} we give an overview of the basic principles of multigrid methods and tensor formats, respectively.  
	Section~\ref{sec:method} gives a detailed explanation on how the individual ingredients of a multigrid method are computed in such a way that the compact format of \eqref{eq:A} is kept. Section~\ref{sec:numtest} includes a variety of numerical tests and a discussion of the efficiency of our multigrid approach for different choices of building blocks and parameters. Concluding remarks and topics for future research are provided in section~\ref{sec:conclusion}.      
	
	\section{Multigrid basics}\label{sec:multigridbasics}
	The building blocks of multigrid methods are
	smoothing schemes, the computation of the set of coarse variables, transfer operators and coarse grid operators. In the following there will be only an overview about what these concepts are and how they work together in a multigrid ansatz. For a detailed treatment and motivation we refer the reader to~\cite{TrottenbergOsterleeSchueller2001}. A prototype of a $V$-cycle multigrid method is given in Algorithm \ref{alg:vZyk}.
	
	\begin{algorithm}
		\DontPrintSemicolon
		$ v_l = \textnormal{MG}(b_l,v_l) $\;
		\uIf{ coarsest grid is reached}{solve coarse grid equation $A_lv_l=b_l$.} 
		\Else{Perform $\nu_1$ smoothing steps for $A_lv_l = b_l$ with initial guess $ v_l $\; 
			Compute the residual $r_l=b_l-A_lv_l$ \;
			Restrict $b_{l+1}=Q_lr_{l}$\;
			$ v_{l+1}=0$\;
			$ e_{l+1}=\textnormal{MG}(b_{l+1}, v_{l+1})$\;
			Interpolate $e_l=P_le_{l+1}$ \;
			$v_l=v_l+e_l$ \;
			Perform $\nu_2$ smoothing steps for $A_lv_l = b_l$ with initial guess $ v_l $\;
		}
		\caption{Multigrid $V$-cycle\label{alg:vZyk}}
	\end{algorithm}
	
	The smoothing process typically consists of a few iterations of a simple iterative method like weighted Jacobi, Gau{\ss}-Seidel or a Krylov subspace method like Richardson or GMRES; cf.~\cite{Greenbaum1997,Hackbusch1994,TrottenbergOsterleeSchueller2001}. Afterwards the error $e$ of the current iterate $v$ is calculated by approximately solving the residual equation $Ae = r$. This is done by performing the computations on a problem of smaller size on a so called coarse grid. To do so, the residual $r$ and the operator $A$ are restricted to a smaller space. In case it is the coarsest level in the hierarchy the restricted system is solved exactly, otherwise the process of smoothing and restriction is repeated until a grid is reached which is small enough for direct computations. 
	Once the coarsest system is solved, the calculated error is interpolated to the next finer grid and added to the current iterate. A final smoothing operation is applied and the process is repeated until the finest grid is reached. 
	
	For notational simplicity we use a two grid notation whenever applicable, i.e., all quantities related to the coarse grid have subscript $c$. In cases where it is necessary to distinguish more than two subsequent grids we use a numbering of the grids as implied by algorithm~\ref{alg:vZyk}.
	
	In order to define restriction and interpolation operators one has to specify coarse variables. 
	There are many different approaches available to define them, e.g., geometric coarsening~\cite{Brandt1977}, compatible relaxation~\cite{Brandt2000, BrannickFalgout2010}, aggregation~\cite{BrezinaManteuffelMcCormichRugeSanders2010} and some others which can be found in~\cite{TrottenbergOsterleeSchueller2001}. Typically these approaches split the variables into two sets, a set of fine variables $\mathcal{F}$ and a set of coarse variables $\mathcal{C}$, which are used in the definition of restriction and interpolation. 
	Assuming that such a splitting is chosen, restriction and interpolation operators are described by rectangular matrices. The restriction operator $Q$ maps the residual to the coarse variables and the interpolation operator maps the calculated error correction from the coarse variables to the original space
	\begin{equation*}
		Q: \mathbb{R}^{|\mathcal{C}\cup\mathcal{F}|} \rightarrow \mathbb{R}^{|\mathcal{C}|}\quad P: \mathbb{R}^{|\mathcal{C}|}\rightarrow\mathbb{R}^{|\mathcal{C}\cup\mathcal{F}|}.
	\end{equation*}
	There are also many different approaches to build these operators, e.g., least squares interpolation~\cite{BrandtBrannickKahlLivshits2011,ManteuffelMcCormickParkRuge2010}, linear interpolation and  others~\cite{TrottenbergOsterleeSchueller2001}. The coarse grid operator $A_c$ is then formed by the Petrov--Galerkin construction $QAP \in \mathbb{R}^{|\mathcal{C}|\times|\mathcal{C}|}$ and the coarse residual is given by $r_c = Qr$. 
	This is the standard choice in algebraic multigrid methods and gives us the approximated residual equation $A_ce_c=r_c$. 
	
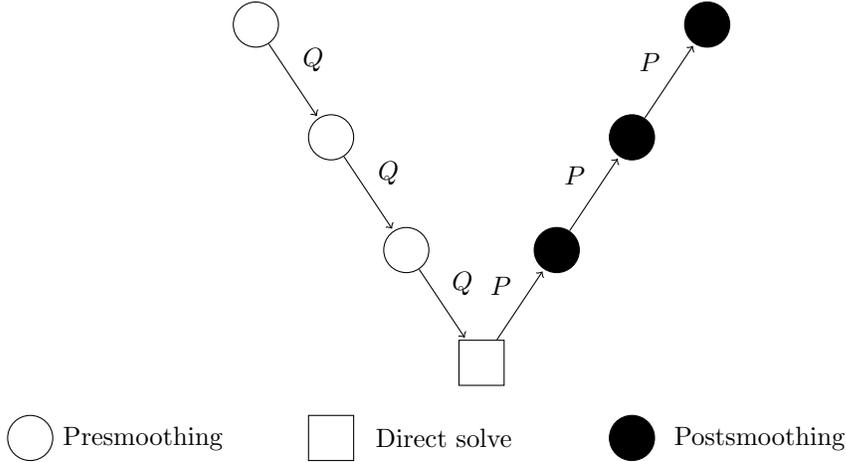
\begin{figure}
		\centering
		\begin{tikzpicture}[shorten >=1pt,auto]
			\draw (2,4.5) node[circle,draw,minimum size = 0.6cm] (A)  {};
			\draw (3,3) node[circle,draw,minimum size = 0.6cm] (B)  {};
			\draw (4,1.5) node[circle,draw,minimum size = 0.6cm] (C)  {};
			\draw (5,0) node[rectangle,draw,minimum size = 0.6cm] (D)  {};
			\draw (6,1.5) node[circle,draw,fill=black,minimum size = 0.6cm] (E) {};
			\draw (7,3) node[circle,draw,fill=black,minimum size = 0.6cm] (F)  {};
			\draw (8,4.5) node[circle,draw,fill=black,minimum size = 0.6cm] (G)  {};
			
			\draw (-1,-1) node[circle,draw,minimum size = 0.6cm] (H)  {};
			\draw (0.5,-1) node (H2)  {Presmoothing}; 
			\draw (3,-1) node[rectangle,draw,minimum size = 0.6cm] (I)  {}; 
			\draw (4.5,-1) node (H2)  {Direct solve};
			\draw (7,-1) node[circle,draw,fill=black,minimum size = 0.6cm] (J)  {};
			\draw (8.7,-1) node (J2)  {Postsmoothing};
			
			\path[->] (A) edge node {$Q$} (B);
			\path[->] (B) edge node {$Q$} (C);
			\path[->] (C) edge node {$Q$} (D);
			\path[->] (D) edge node {$P$} (E);
			\path[->] (E) edge node {$P$} (F);
			\path[->] (F) edge node {$P$} (G);
			
		\end{tikzpicture}
		\caption{Multigrid V-cycle: On each level, a presmoothing iteration is performed, then the problem is restricted to the next coarser grid. On the smallest grid, the problem is solved exactly by a direct solver. When interpolating back to the finer grids, postsmoothing iterations are applied on each level.}
		\label{fig:vcycle}
	\end{figure}

	Instead of stopping on the second grid and solving $A_ce_c = r_c$ exactly, e.g., because $A_c$ is still to large, one can again solve this system of equations by a two-grid approach. Iterating this idea ultimately yields a multigrid method, where 
	only on the coarsest grid, i.e., the one with the smallest dimension, the corresponding system is solved exactly. 
	This strategy is described in algorithm \ref{alg:vZyk} and gives rise to a $V$-cycle depicted in figure~\ref{fig:vcycle}. Other cycling strategies 
	like $W$- or $F$-cycles~\cite{TrottenbergOsterleeSchueller2001} are also possible, but we do not consider them for the sake of simplicity. However, all ideas developed in this paper can also be used with these cycling strategies in a straight forward way.
	
	\section{Tensor basics}\label{sec:tensorbasics}
	We define a tensor $\mathcal{X}$ as follows.
	\begin{definition} \label{def:tensor}
		Consider the space $\mathbf{R}=\mathbb{R}^{n_1} \otimes \cdots \otimes \mathbb{R}^{n_J}$ which is spanned by 
		\begin{equation}\label{eq:tensorspan}
			\{v^{(1)} \otimes v^{(2)} \otimes \cdots \otimes v^{(J)}: v^{(j)} \in \mathbb{R}^{n_j}, 1 \leq j \leq J \}.
		\end{equation}
		Then each element $\mathcal{X}  \in \mathbf{R}$ is a $J$-way tensor and the generating products $v^{(1)} \otimes v^{(2)} \otimes \cdots \otimes v^{(J)}$ are called elementary tensors. Note, that one can interpret $\mathcal{X}$ as a multidimensional array. We therefore label its entries by $J$ indices as $\mathcal{X}(i_1,\dots,i_J)$.
		
		The rank of a Tensor $\mathcal{X}$ is the smallest integer $r$ that satisfies
		\begin{equation}\label{eq:cp}
			\mathcal{X}=\sum\limits_{i=1}^r \mathcal{V}^{(i)}
		\end{equation}
		for elementary tensors $\mathcal{V}^{(i)}$. 
	\end{definition} 
	
	Note, that a one-way tensor is simply a vector and a two-way tensor is a matrix. Any elementary tensor has tensor rank 1.
	
	We denote the mapping that maps a tensor $\mathcal{X}$ to the corresponding vector by $\tenvec(\mathcal{X})$. The ordering in which the entries appear in the vector is not crucial as long as it is consistent if there are multiple occurrences.   
	
	Closely related to tensors is the concept of the Kronecker product. As the Kronecker product of $J$ vectors $v^{(1)},v^{(2)}, \dots ,v^{(J)}$, a vector of dimension $n_1n_2 \dots n_J$, and the corresponding tensor $v^{(1)} \otimes v^{(2)} \otimes \cdots \otimes v^{(J)}$ have the same entries we use the same symbol $\otimes$ for tensor and Kronecker products. A useful property of the Kronecker product is given by
	\begin{equation}\label{eq:kron}
		(A \otimes B)(C \otimes D) = AC \otimes BD. 
	\end{equation} 
	when $A,B,C,D$ are matrices such that all matrix-matrix products in \eqref{eq:kron} are defined.
	
	The relation~\eqref{eq:kron} is of particular interest in the situation that an elementary tensor 
	$\mathcal{X}=x^{(1)} \otimes x^{(2)} \otimes \cdots \otimes x^{(J)}$ and a Kronecker product of matrices $A^{(1)} \otimes \cdots \otimes A^{(J)}$ denoted by $A$ are given such that $A^{(j)}x^{(j)}$ is well-defined. In this situation we find for the multiplication $A\mathcal{X}$,
	\begin{equation}
		A\mathcal{X} = A^{(1)} x^{(1)} \otimes A^{(2)} x^{(2)} \otimes \cdots \otimes A^{(J)} x^{(J)}.
	\end{equation}
	Thus the product can again be described by an elementary tensor and it can be computed very efficiently by only computing matrix-vector products of small size $n_j$.
	
	Motivated by the simple matrix-tensor multiplication for elementary tensors it is clear that it is beneficial to work with a tensor representation of $\mathcal{X}$ instead of $\tenvec(\mathcal{X})$, i.e., the full representation. There are many tensor formats (e.g., CP~\cite{CarrollChang1970,Kiers2000,KoldaBader2009}, H-Tucker~\cite{Grasedyck2010,HackbuschKuehn2009}, Tensor Train~\cite{Oseledets2010,Oseledets2011,OseledetsDolgov2012}) with different approximation properties and matrix-tensor multiplications. Of these formats we focus on the Tensor Train (TT) format in the following and introduce its basic concepts next. 
	
	\subsection{Tensor Train}\label{subsec:tensortrain}
	A Tensor $\mathcal{X}$ is in TT-format with TT-ranks $r_0,\dots, r_J$ if each entry of $\mathcal{X}$ is given as
	\begin{equation}\label{eq:tt}
		\mathcal{X}(i_1,\dots, i_J) = G_1(i_1)\cdot G_2(i_2)\cdot{}\dots{}\cdot G_J(i_J)
	\end{equation} 
	with parameter dependent matrices $G_k(i_k) \in \mathbb{R}^{r_{k-1} \times r_k}, k=1, \dots, J$. 
	As $\mathcal{X}(i_1,\dots, i_J)$ is a scalar value, \eqref{eq:tt} implies that $r_0 = r_J=1$.
	Note that the $G_k$ can be interpreted as $r_{k-1} \times n_k \times r_k$ tensors and they are called cores of $\mathcal{X}$. 
	The storage complexity of $\mathcal{X}$ in format \eqref{eq:tt} is bounded by $(J-2)\widehat{n}\widehat{r}^2+2\widehat{n}\widehat{r}$, if $n_k \leq \widehat{n}$ and $r_k\leq \widehat{r}$ for $k=1, \dots , J$.
	
	To be able to use the TT-format in an iterative method, basic vector operations should perform cost efficiently. 
	The addition of two tensors in the TT-format can be computed as follows:
	Given tensors $\mathcal{X}$ and $\widetilde{\mathcal{X}}$ with cores $G_k(i_k), \widetilde{G}_k(i_k)$ respectively, then the cores of $\mathcal{X} + \widetilde{\mathcal{X}}$ are given by
	\begin{eqnarray*}
		\widehat{G}_1(i_1) &=& \left( G_1(i_1)\quad \widetilde{G}_1(i_1) \right),\qquad\widehat{G}_J(i_J) = \left(\begin{array}{c} G_J(i_J) \\ \widetilde{G}_J(i_J)\end{array}\right), \\
		\widehat{G}_k(i_k) &=& \left(\begin{array}{cc} G_k(i_k) & 0 \\ 0 & \widetilde{G}_k(i_k)\end{array}\right) \textnormal{ for } k = 2,\dots, J-1.
	\end{eqnarray*}
	Thus, no arithmetic operations are necessary, but the TT-ranks (and therefore the storage complexity) of $\mathcal{X}+\widetilde{\mathcal{X}}$ are higher than those of $\mathcal{X}$ and $\widetilde{\mathcal{X}}$. Therefore, a rounding (or truncation) procedure is necessary, which allows to approximate a given TT-tensor by a TT-tensor of lower rank. This truncation can be done in a numerically stable way with cost $\mathcal{O}(J\widehat{n}\widehat{r}^3)$ by using~\cite[Algorithm~2]{Oseledets2010}. Another operation which is typically needed in iterative methods is the scaling of a vector (or tensor) by a scalar $\alpha$. This can be done efficiently (and without changing the TT-ranks) by multiplying one of the cores by $\alpha$.
	Another operation needed in our method is the Euclidean inner product. 
	It can be evaluated by~\cite[Algorithm 4]{Oseledets2010} with cost $\mathcal{O}(J\widehat{n}\widehat{r}^3)$.
	
	A matrix $A$ is in TT-format if each entry can be computed as
	\begin{equation}\label{eq:tt_matrix}
		A(i_1,\dots, i_J; j_1,\dots, j_J) = M_1(i_1,j_1)\cdot M_2(i_2,j_2)\cdot{}\dots{}\cdot M_J(i_J, j_J)
	\end{equation}
	with parameter dependent matrices $M_k(i_k,j_k) \in \mathbb{R}^{r_{k-1} \times r_k}, k = 1,\dots, J$. Note that a matrix $A$ with TT-ranks $r_0 = r_1 = \dots r_J = 1$ is just the Kronecker product of $J$ matrices, i.e.,
	\begin{equation}
		A = M_1 \otimes M_2 \otimes \dots \otimes M_J.
	\end{equation} 
	Given a matrix $A$ and a tensor $\mathcal{X}$ in TT-format with all TT-ranks bounded by $r$, the matrix-vector product can be computed with arithmetic cost $\mathcal{O}(J\widehat{n}^2\widehat{r}^4)$ by using~\cite[Algorithm 5]{Oseledets2010}. As the TT-ranks of $A$ and $\mathcal{X}$ multiply when performing this operation, it is crucial to use truncation afterwards.
	
	In \cite[Section 3]{Oseledets2010} it is described how to efficiently convert a tensor in CP decomposition into the TT-format. Similar techniques can be used to convert a matrix $A$ in representation \eqref{eq:A} into the TT-format.
	
	\section{Method}\label{sec:method}
	Tensor-based methods are efficient only when the tensors appearing in the problem have low rank. If this is not the case, the tensor structure \eqref{eq:A} of $A$ can be used to implement fast matrix vector multiplication, but additional savings by using one of the described tensor formats for the iterate are not possible. Thus the success of our attempt to use a tensor format for the iterates depends crucially on the rank of the steady state vector. 
	For non-interacting systems the stationary distribution has rank one, as stated by the following theorem, which is a basic result from linear algebra.
	
	\begin{theorem}\label{the:noninteracting}
		For a non-interacting system, i.e., $E_j^t=I$ for $t \neq j$ and $T=J$ in \eqref{eq:A} the stationary distribution $p$ satisfies 
		\begin{equation}
			p=\bigotimes\limits_{j=1}^J p^{(j)},
		\end{equation}
		where $p^{(j)}$ is the stationary distribution of $E_j^j$, i.e., it corresponds to a tensor of rank one.
	\end{theorem}
	
	The proof is simple, as the tensor product of the individual stationary distributions of the systems is the stationary distribution of the combined system.
	
	Theorem~\ref{the:noninteracting} suggests that for systems with a low level of interaction the stationary distribution can be approximated by a tensor with low rank. Numerical experiments in, e.g.,~\cite{KressnerMacedo2014} as well as our experiments in section \ref{sec:numtest} show that in practice a rather low rank is often sufficient. Thus savings of both memory and compute time can be realized by using a truncated tensor format in an iterative method.
	
	This result as well as the numerical experiments motivate the development of multigrid methods for system matrices $A$ with tensor structure that allow for the efficient use of tensor formats as described in section~\ref{sec:tensorbasics}. In the following we will specify the building blocks of a multigrid method according to the description in section~\ref{sec:multigridbasics}. Our main goal in the following is to preserve the tensor structure of $A$ in~\eqref{eq:A}.
	
	\subsection{Smoother}
	Unlike in multigrid methods for matrices the tensor structure of the linear systems at hand rules out smoothers that require access to individual entries of the operator, e.g., Jacobi, which requires the diagonal, or Gau{\ss}-Seidel, which requires the lower triangular part. Unless there is a way to represent these parts of the operator in a tensor format as well, it is very difficult and prohibitively expensive to use these smoothers. Thus typical candidates that can be applied in a tensor environment are smoothers that only require multiplications of the iterate with the system matrix, scalar products, multiplication with scalars and addition of vectors (or their tensor representations). These requirements are fulfilled by, e.g., Richardson method or more generally any polynomial method and thus also by Krylov subspace methods. As the optimal smoothing parameter in the case of Richardson is hard to find we suggest to use a Krylov subspace method like GMRES with a fixed number of iterations as the prototype smoother in a tensor environment. As long as the number of iterations done in the smoother is small the overhead in storage and computation required for GMRES are negligible.
	
	In case the system matrix allows for an explicit representation of the diagonal $D$ or the upper and lower triangular matrices $L$ and $U$ with zero diagonal such that $A = D - L -U$, a Jacobi or Gau{\ss}-Seidel like smoother can be employed by approximately inverting $(D - L)$ using a method like AMEn~\cite{DolgovSavostyanov2013a,DolgovSavostyanov2013b}. The effectiveness of such an approach hinges mainly on the conditioning of the linear systems that need to be solved in AMEn. Especially as it is known that the systems arising in AMEn for $A$ become increasingly ill-conditioned for growing problem sizes, special care has to be paid when applying this approach.
	
	\subsection{Selection of coarse variables}
	The following proposition gives us a starting point for achieving our goal of maintaining the structure of $A$. For that purpose we assume that the restriction operator $Q$ and the prolongation operator $P$ itself possess the same tensor structure as $A$. The proposition can be proven using basic properties of the Kronecker product.
	\begin{proposition}\label{pro:QAP}
		Let $A$ of the form \eqref{eq:A} be given with $E_j^t \in \mathbb{R}^{n_j \times n_j}$. Let $P = \bigotimes_{j = 1}^J P_j$ and $Q = \bigotimes_{j = 1}^J Q_j$ with $P_j \in \mathbb{R}^{n_j \times n_j^c}$ and $Q_j \in \mathbb{R}^{n_j^c \times n_j}$ where $n_j^c << n_j$. Then the corresponding Petrov--Galerkin operator satisfies
		\begin{equation}\label{eq:QAP}
			QAP = \sum\limits_{t = 1}^T\bigotimes_{j = 1}^J Q_j E_j^t P_j.
		\end{equation}
	\end{proposition}
	
	A direct consequence of this proposition is that the set of coarse variables can be computed for the smaller problems described by the matrices $E_j^{\star}$. However, it is not always straightforward how to do this. In case the matrices $E_j^{\star}$ originate from some geometric structure, like the overflow queuing network example in section \ref{subsec:overflow}, one can use an established coarsening algorithm. A matrix like in the Kanban system example in section \ref{subsec:kanban} does not correspond to a connected graph, indeed there are several isolated states. In this case, coarsening can be done with respect to the auxiliary matrices 
	\begin{equation}\label{eq:E_tilde}
		\widetilde{E}_j = \sum\limits_{t = 1}^T E_j^t,
	\end{equation}
	which collect all transitions of the $j$-th subsystem and correspond to connected graphs again.
	Problems of the first kind can be thought of $j$ individual Markov chains which are coupled in some way, but are also meaningful on their own, while the second kind are systems which have a high level of interaction where each individual subsystem is dysfunctional without interacting with the other systems. This will be seen in more detail in the description of the different test examples in section \ref{sec:numtest}.
	
	\subsection{Transfer operators}
	Once the set of coarse level variables is selected, transfer operators $P$ and $Q$ have to be defined in accordance with our assumption on their structure.
	The way of forming $P$ and $Q$ as described in Proposition~\ref{pro:QAP} directly implies that they are built with respect to the ``local'' operators $E_j^{\star}$ of smaller dimension and lifted to the full dimension of $A$ via Kronecker products.
	The choice of the ``local'' interpolation and restriction operators depends on the considered problem: In case the problem is similar to the discretization of a partial differential equation, linear interpolation and a $J$-fold Kronecker product thereof is a reasonable choice. This is the case for the overflow queuing network problem considered in section~\ref{sec:numtest} which bears similarity to an advection-diffusion problem. 
	In principle all known constructions from algebraic multigrid can be used to construct the transfer operators, including adaptive techniques.	
	
	\subsection{Coarse grid operator}
	The coarse grid operator is chosen as the Petrov--Galerkin operator, i.e., $A_c = QAP$, being the natural extension of the Galerkin operator that is the optimal coarse grid operator in the symmetric case.
	
	Finally, on the coarsest grid the residual equation has to be solved. There are two ways to do this, either the coarsening results in a problem of a size where direct solving without exploiting the tensor structure is possible or by an iterative method (e.g., GMRES) that exploits the tensor structure. In both cases truncation has to be applied to obtain a low rank tensor representation of the solution.
	
	In summary, by keeping in mind these considerations, we are able to maintain the tensor structure of the system matrix $A$ throughout all levels of our multigrid method. This structure has two main advantages: First, only the small matrices need to be stored, substantially reducing the storage requirements and second an efficient matrix vector product can be performed if $\mathcal{X}$ is given in a tensor format like the TT-format described in section \ref{subsec:tensortrain} with low rank. While the chosen interpolation and restriction operators are of rank 1 the system matrix A usually is not. The first implies that multiplication with the grid transfer operators does not change the rank of a vector. The same holds for the coarse grid operator that is formed as the Petrov--Galerkin product. Though multiplication with the system matrix or the coarse grid operators usually results in a larger rank, as the resulting rank of a matrix vector product $A\mathcal{X}$, where both $A$ and $\mathcal{X}$ are in TT-format with ranks $r_A$ and $r_\mathcal{X}$ respectively, is the product of the ranks $r_A r_\mathcal{X}$. Therefore, after each multiplication with an operator of rank larger than 1 truncation is employed.
	
	\section{Numerical Tests} \label{sec:numtest}
	In this section we illustrate how to choose the ingredients of a multigrid approach based on two examples. All computations were performed in MATLAB R2013a using the TT-Toolbox~\cite{Toolbox}.
	
	\subsection{Overflow queuing network}\label{subsec:overflow}
	The first example we consider is the so called Overflow Queuing Network.
	A finite number $J$ of queues is given with a corresponding capacity $k_i,\, i = 1,\dots,J$.
	Customers arrive to an arbitrary queue $q_i,\, i = 1,\dots,J$ according to a Poisson process with rate $\lambda_i,\, i = 1,\dots,J$ and are served with an exponentially distributed service time with rate $\mu_i,\, i = 1,\dots,J$. These local processes can be described by the following matrices
	\begin{equation}\label{eq:oqn_subproblem}
		E_i^{(i,i)} = \begin{pmatrix}
			0&\mu_i&&&0\\
			\lambda_i& \ddots & \ddots\\
			&\ddots&\ddots&\ddots&\\
			&&\ddots&\ddots&\mu_i&\\
			0&&&\lambda_i&0\\
		\end{pmatrix},\, i = 1,\dots, J.
	\end{equation}
	Note that here and in the following, for a better understanding we denote the transitions by two-tuples instead of natural numbers. Here, the tuple $(i,j)$ describes a transition which is initiated by queue $i$ and affects queue $j$. The synchronized events appear if queue $q_i$ has reached its capacity, i.e., $q_i$ is full. In this situation an arriving customer has to enter $q_{i+1}$. If $q_{i+1}$ is also full, the customer will enter the subsequent queue $q_{i+2}$ and so on. If all of the subsequent queues are full, the customer will leave the system. For sake of simplicity we only describe the synchronized events in which the customer can enter the subsequent queue, because it is not full:
	\begin{itemize}
		\item if $q_i$ is full, the customer leaves $q_i$, which leads to the following matrix 
		$$
		E_i^{(i,i+1)} = \begin{pmatrix}
			0&0&&&\\
			0& \ddots & \ddots&&\\
			&\ddots&\ddots&\ddots&&\\
			&&\ddots&\ddots&0&0&\\
			&&&&0&\lambda_i\\
		\end{pmatrix},\, i= 1, \dots,J-1.
		$$
		\item if $q_i$ is full, $q_{i+1}$ gets an additional customer, which leads to the following matrix 
		$$ E_{i+1}^{(i,i+1)} = \begin{pmatrix}
			0&0&&&0\\
			1& \ddots & \ddots\\
			&\ddots&\ddots&\ddots&\\
			&&\ddots&\ddots&0&\\
			0&&&1&0\\
		\end{pmatrix},\, i = 1,\dots, J-1.
		$$
	\end{itemize}
	
\begin{figure}[t]
		\centering
		\subfloat [][]{
			\includegraphics[width=.32\textwidth]{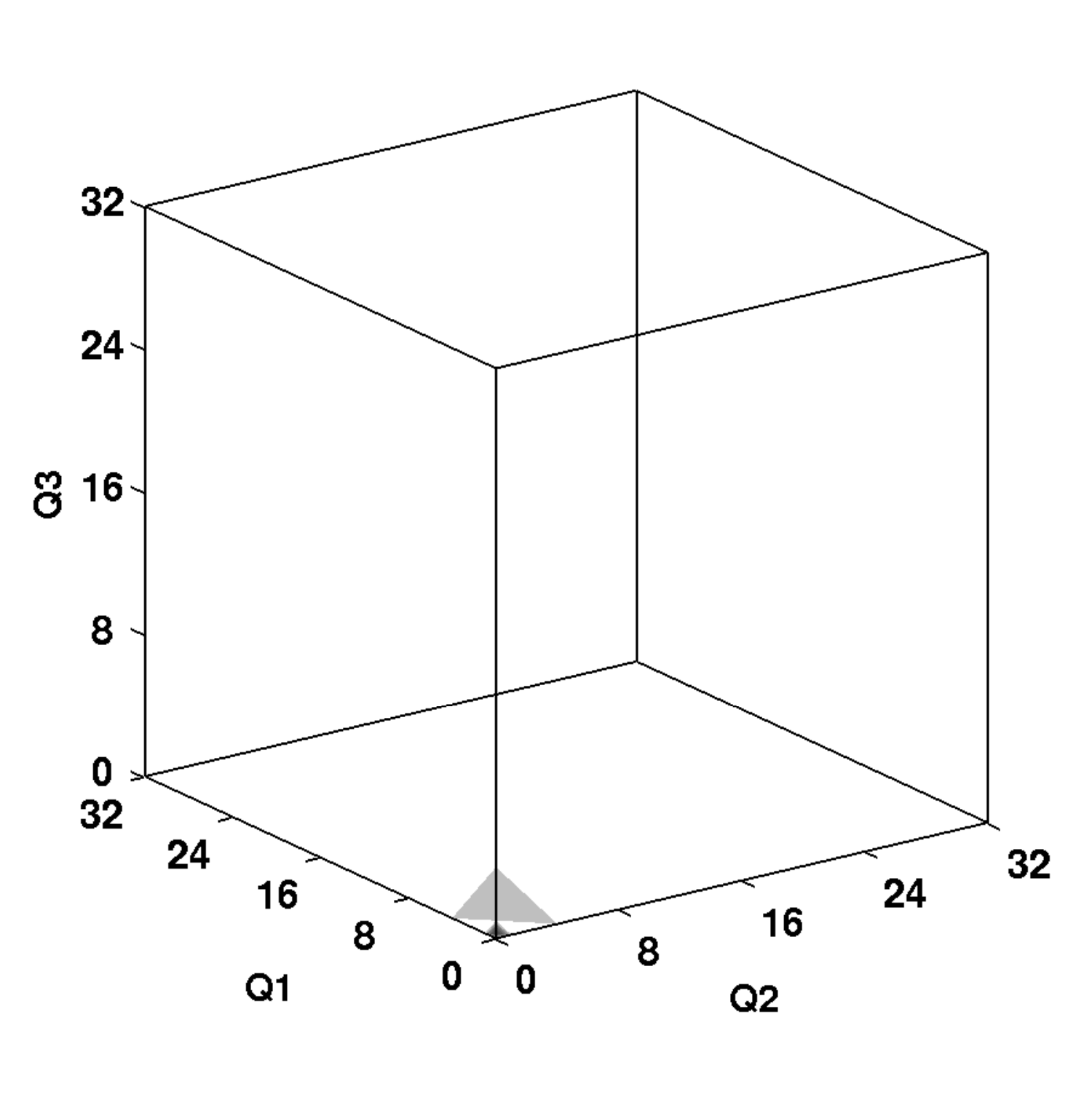}
			\label{fig:solplot1}
		}\hspace{-.04\textwidth}
		\subfloat[][]{
			\includegraphics[width=.32\textwidth]{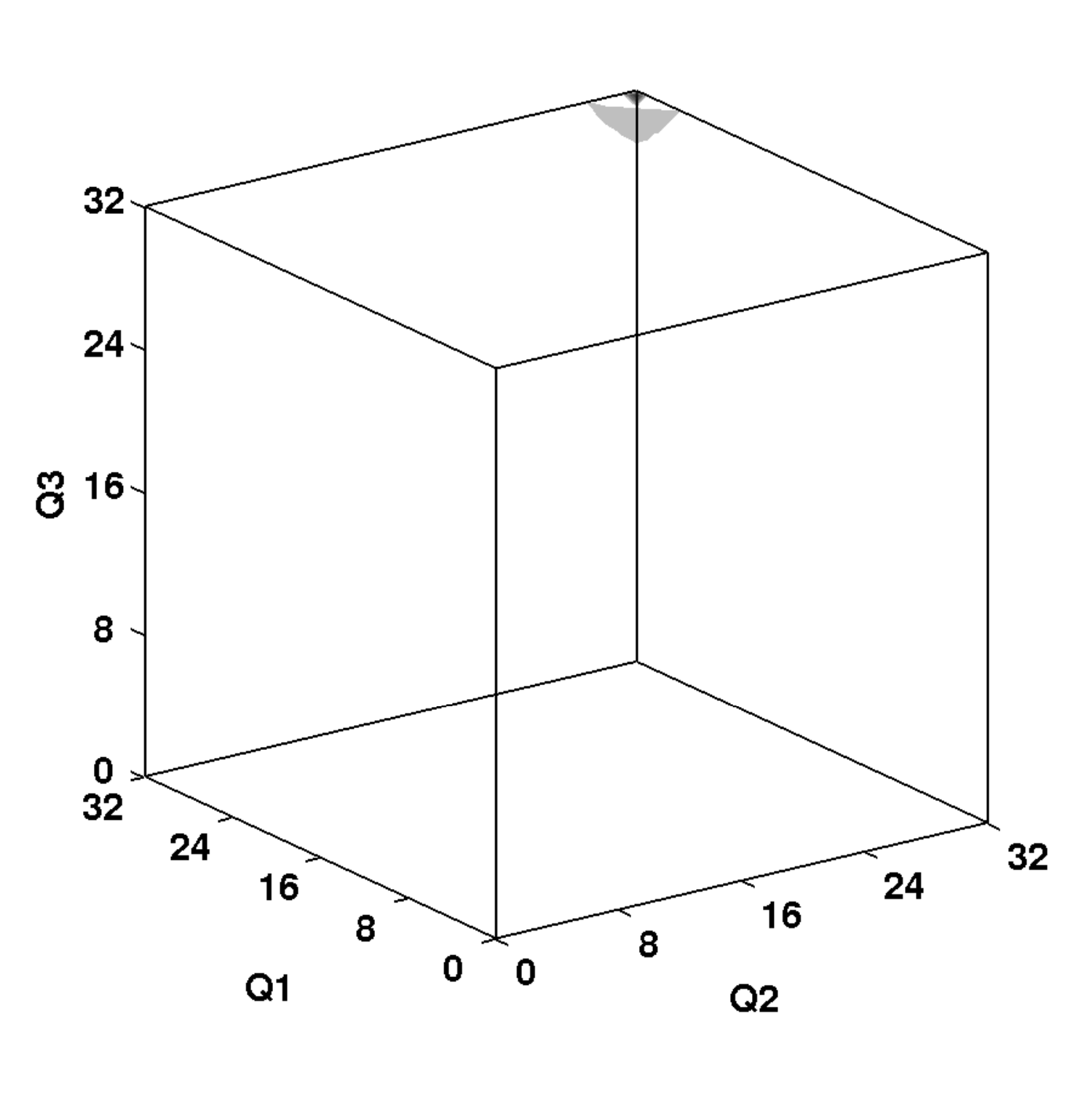}
			\label{fig:solplot2}
		}\hspace{-.04\textwidth}
		\subfloat[][]{
			\includegraphics[width=.32\textwidth]{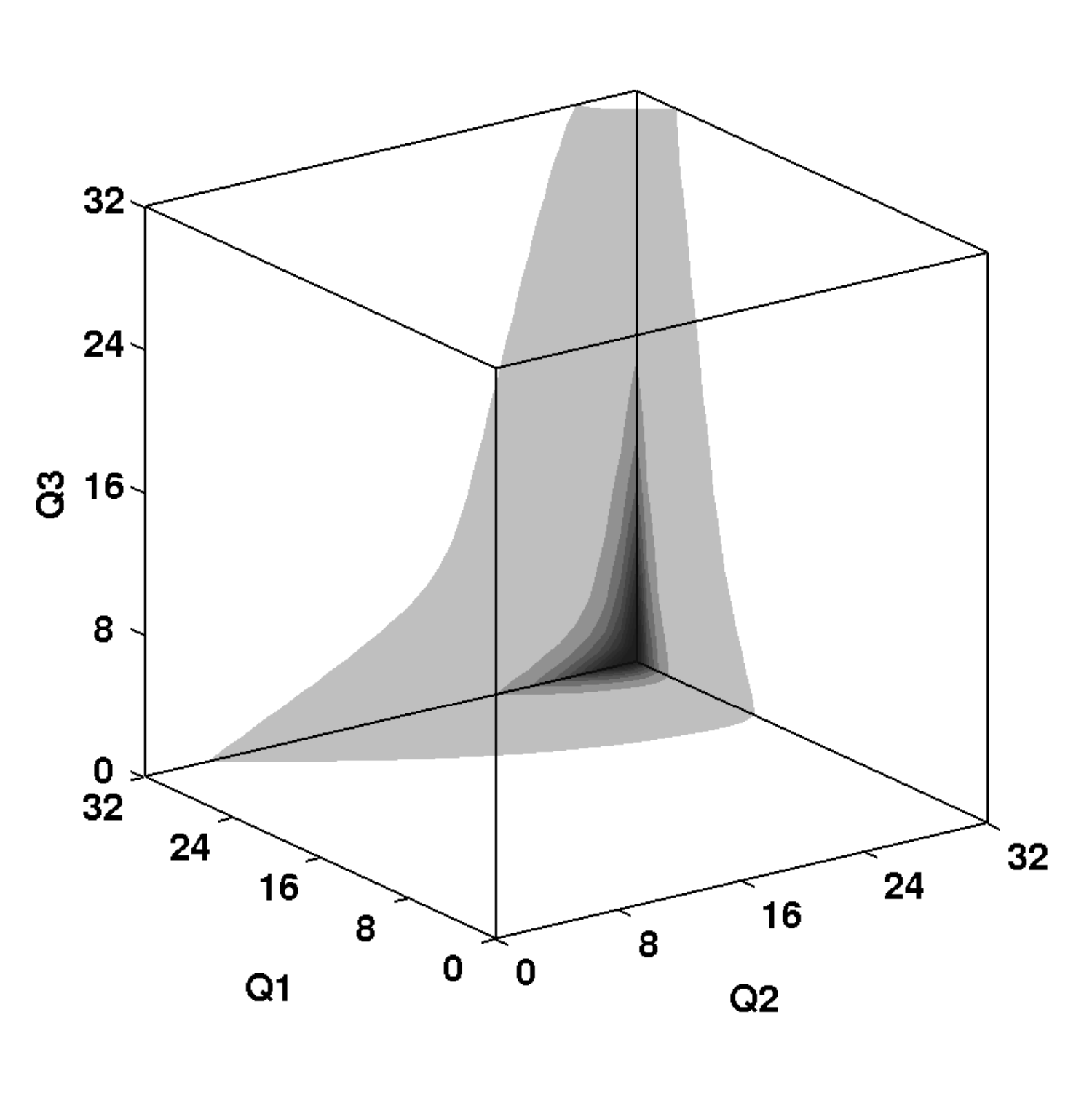}
			\label{fig:solplot3}
		}
		\caption{Solution for parameters (a) $\mu=[1, 1, 1] \text{ and } \lambda=[0.1, 0.1, 0.1]$, (b) $\mu=[0.1, 0.1, 0.1] \text{ and } \lambda=[1, 1, 1]$, (c) $\mu=[0.25, 0.5, 1] \text{ and } \lambda=[0.5, 0.5, 0.5]$.}
		\label{fig:solplot}
	\end{figure}

\begin{figure}[t]
		\centering
		\subfloat [][]{
			\includegraphics[width=.32\textwidth]{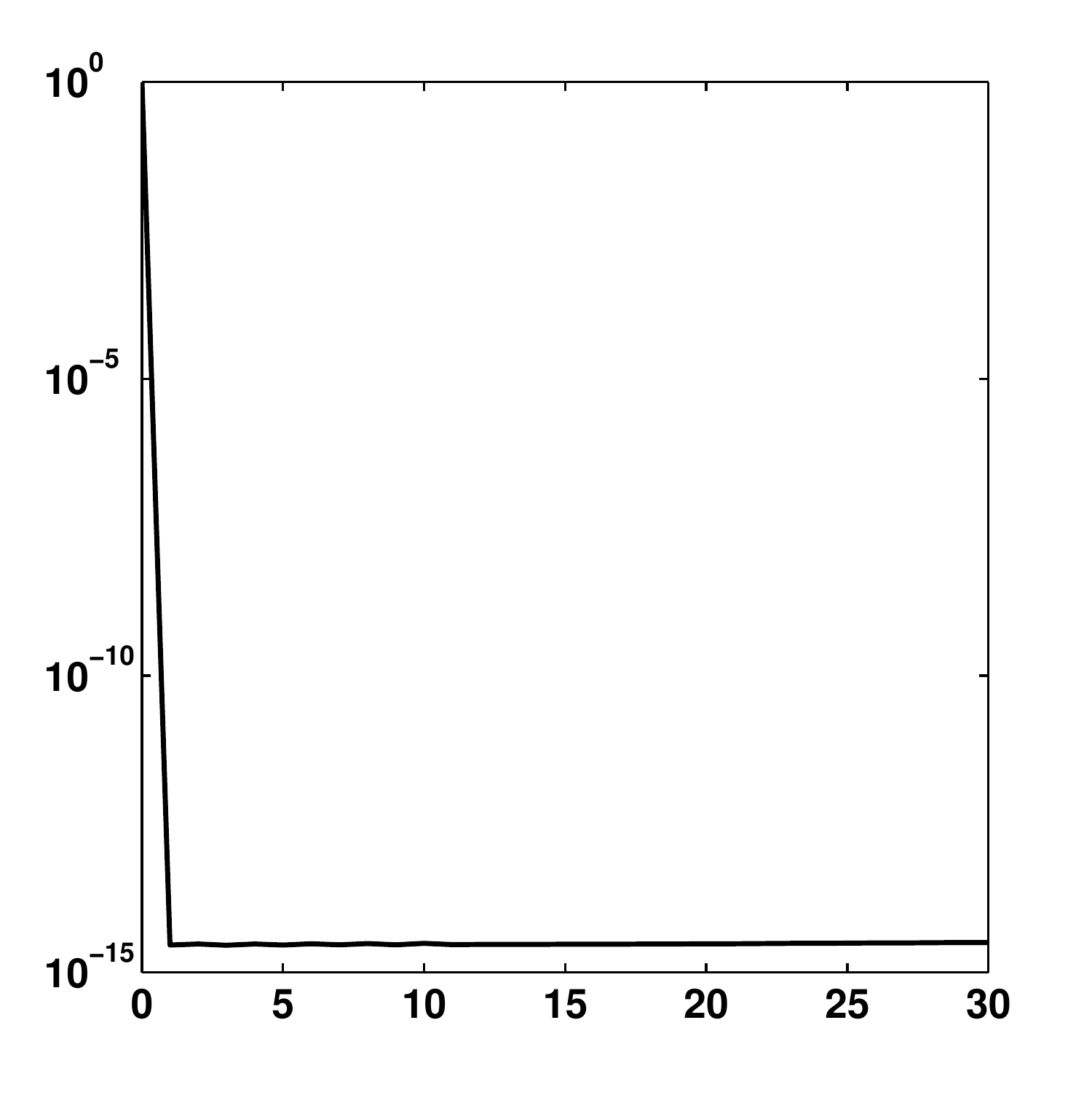}
			\label{fig:rankplot1}
		}
		\subfloat[][]{
			\includegraphics[width=.32\textwidth]{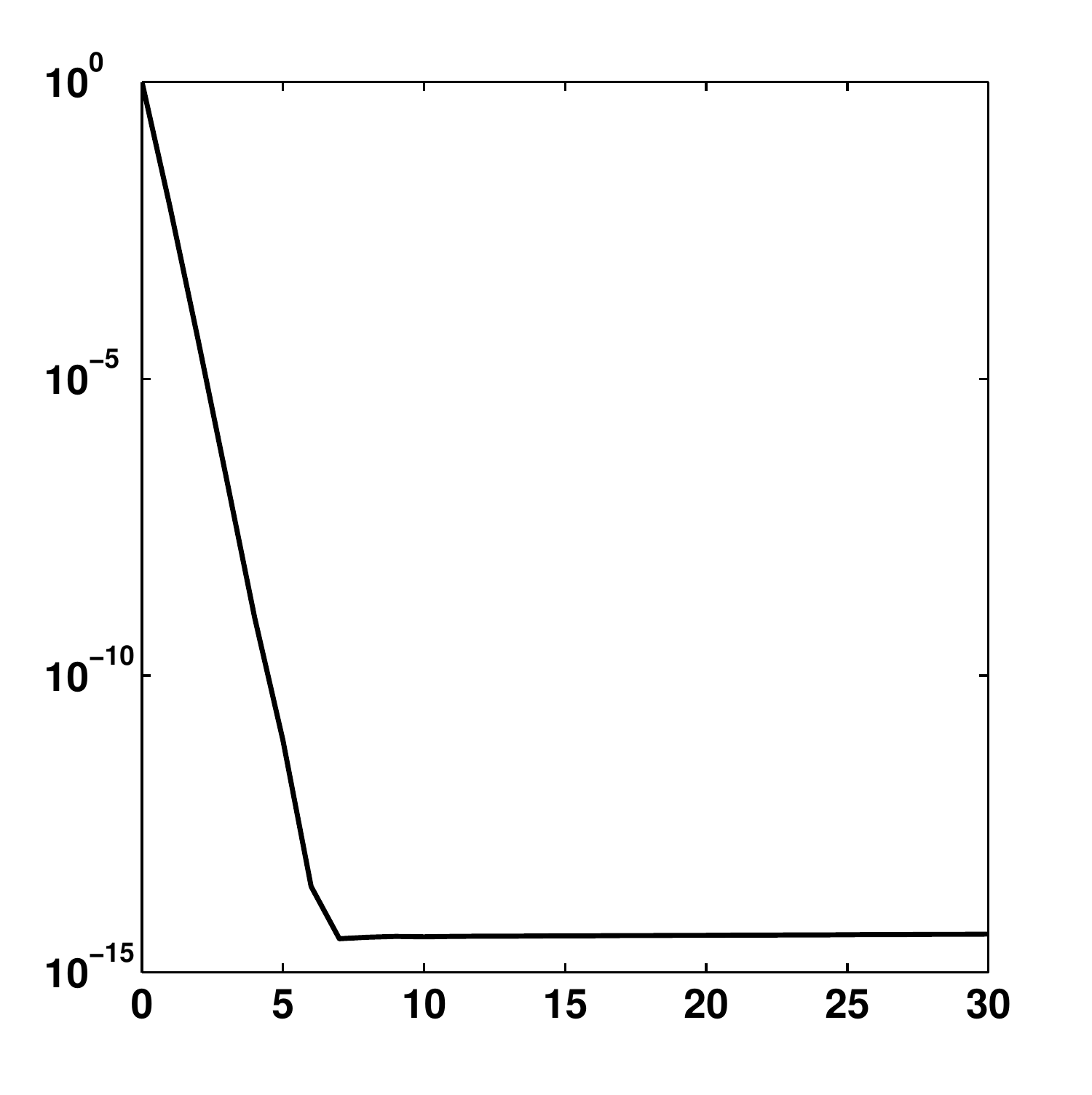}
			\label{fig:rankplot2}
		}
		\subfloat[][]{
			\includegraphics[width=.32\textwidth]{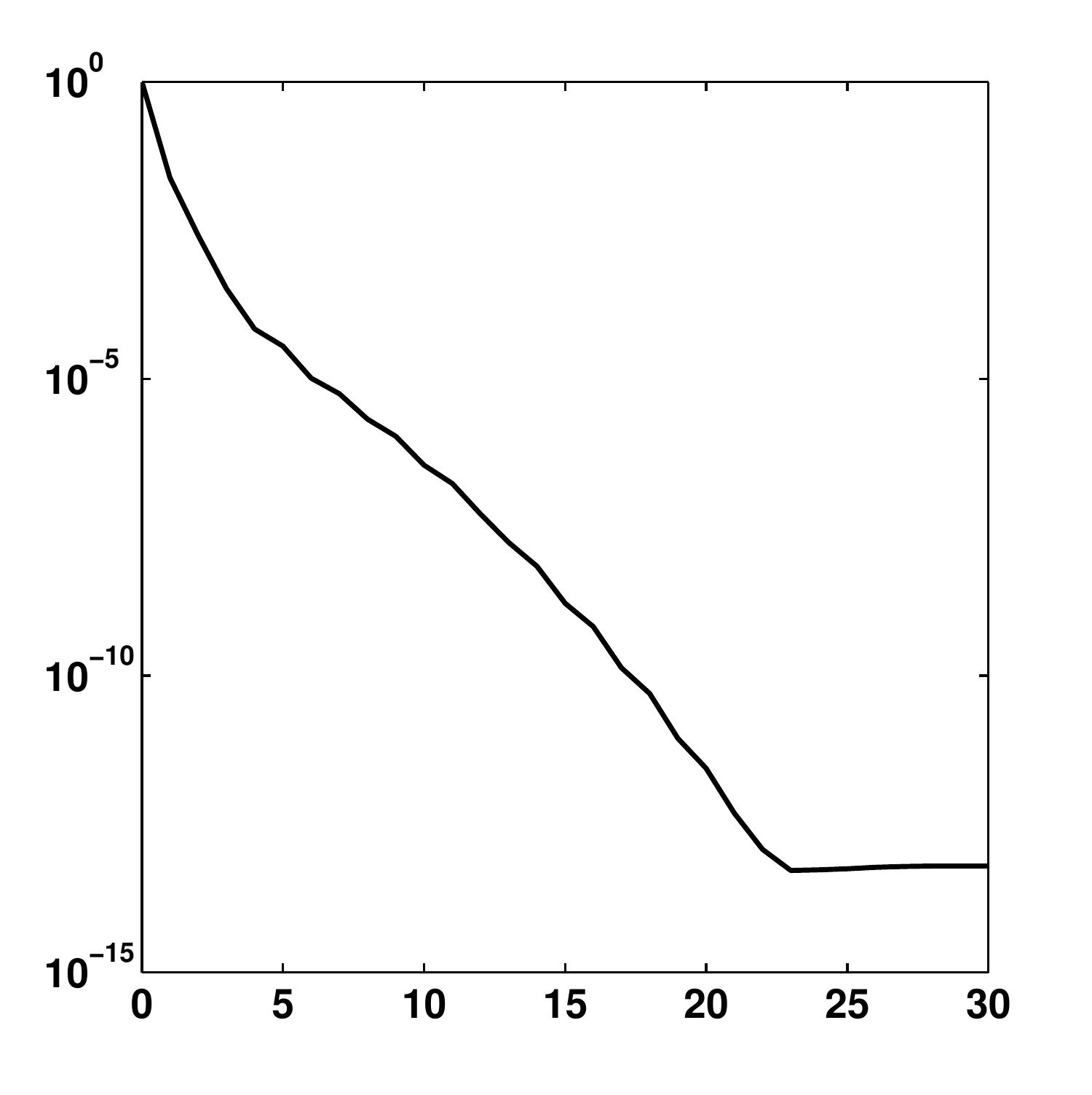}
			\label{fig:rankplot3}
		}
		\caption{Accuracy of rank $R$ approximations (computed by truncation of the exact solution) for parameters (a) $\mu=[1, 1, 1] \text{ and } \lambda=[0.1, 0.1, 0.1]$, (b) $\mu=[0.1, 0.1, 0.1] \text{ and } \lambda=[1, 1, 1]$, (c) $\mu=[0.25, 0.5, 1] \text{ and } \lambda=[0.5, 0.5, 0.5]$.}
		\label{fig:rankplot}
	\end{figure}

	All of the matrices $E_i$ are of the dimension $n_i=k_i+1$. For a detailed description of the model see~\cite{Buchholz2000}. Note that the choices of the parameters $\lambda$ and $\mu$ affect whether more synchronized or local events appear.
	For a better understanding, consider 3 queues with capacity $32$ and different choices of the parameters $\mu=[\mu_1, \dots, \mu_J]$ and $\lambda=[\lambda_1,\dots, \lambda_J]$, which can be seen in figure \ref{fig:solplot}.
	The figures \ref{fig:solplot1}, \ref{fig:solplot2} and \ref{fig:solplot3} show the distributions of the solutions for the corresponding parameter choice on each queue. In every figure we have three axes (one for each queue) and each axis describes the capacity.
	In case customers arrive slowly and are served rapidly (small $\lambda$, large $\mu$) all queues are (almost) empty; cf.~figure~\ref{fig:solplot1}. On the other hand, if customers arrive in short succession and are served slowly we obtain a steady state solution that corresponds to three full queues as can be seen in figure~\ref{fig:solplot2}. Both parameter sets have limited interaction, in the first situation no queue ever flows over which results in no interaction and in the second case almost every customer that is rejected at a queue is rejected at each subsequent queue as well, resulting effectively in no interaction as well.
	Only if the parameters are chosen such that slow and fast queues are mixed the steady state distribution becomes non-trivial as can be seen in~\ref{fig:solplot3}.
	In figure \ref{fig:rankplot} we investigate fixed rank approximations of the solutions for the three different problems and observe that for the solution \ref{fig:solplot3} a higher rank is needed than for the other ones due to the increased interaction of the queues. This observation confirms the assumption that a higher level of interaction in the model leads to a higher rank for its solution. Nevertheless, an accurate approximation of the solution is still obtainable with low rank.
	
%

	\begin{figure}[t]
		\centering
		\includegraphics[width=.6\textwidth]{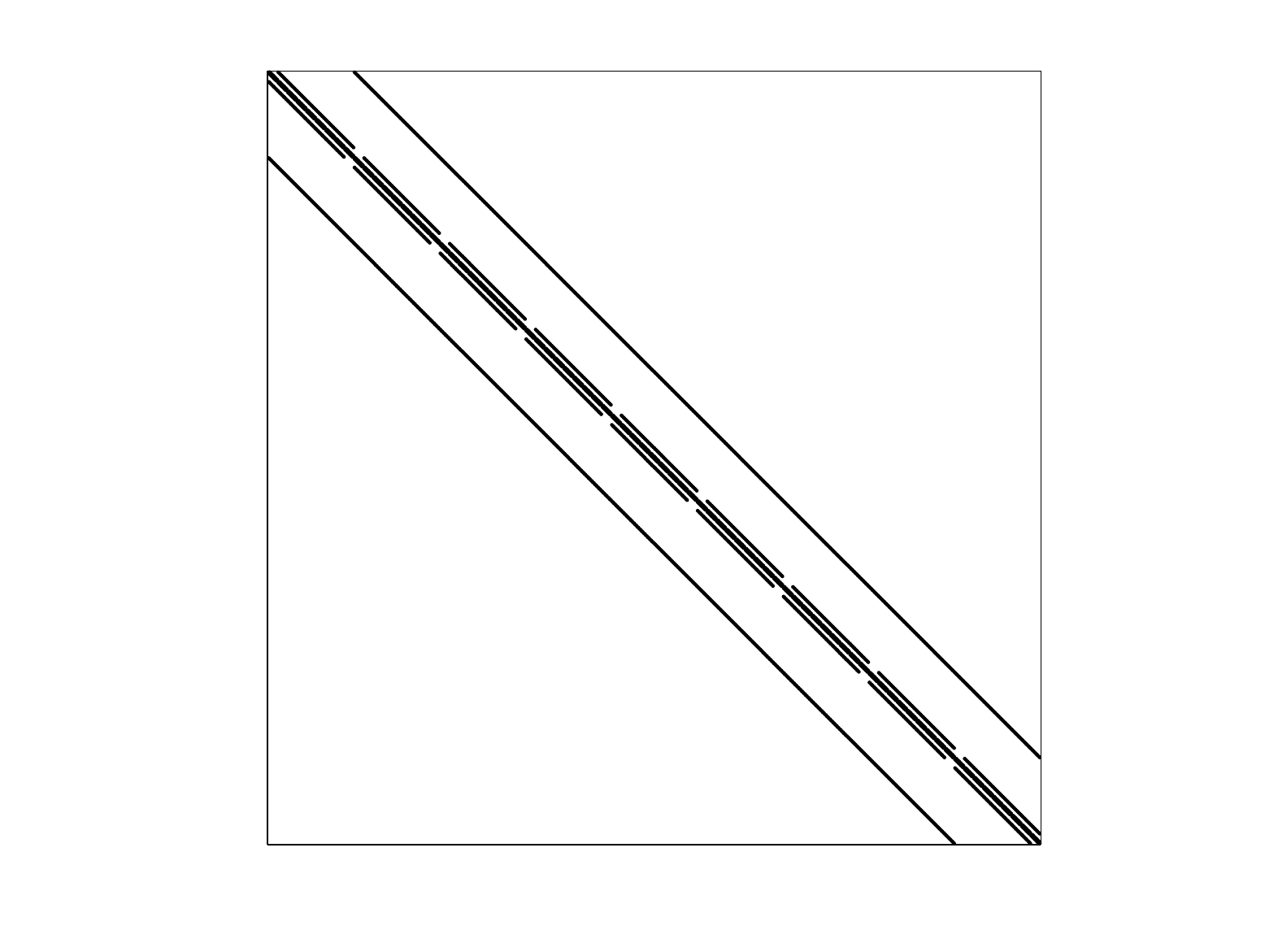}
		\caption{Nonzero structure of overflow queuing network for $J = 3, k = [8, 8, 8]$.}
		\label{fig:spyQ}
	\end{figure}
	
	One way to choose the ingredients for a multigrid approach is to take into account the underlying geometric structure and to preserve some important properties of the problem matrix $A$ on the coarse grids. In our case, the property that our matrix $A$ has column sum zero should be kept. The structure of our overflow queuing problem for $J = 3$, $k_i = 8$, $i = 1,2,3$ is shown in Figure~\ref{fig:spyQ}. Note that the nonzero structure is exactly the same as that of the standard finite difference discretization of the $J$-dimensional Laplace operator.
	A look at the geometric structure of the local problem in an overflow queuing problem~\eqref{eq:oqn_subproblem} and the structural similarity to the Laplace operator motivated us to use full coarsening and direct interpolation for the interpolation operator, based on these local matrices, for each queue. To keep the column sum of the coarse grid matrix zero, we use linear restriction, because our coarse grid matrix is a Petrov--Galerkin operator and we then have $\mathbf{1}^TQAP = \mathbf{1}^TAP = 0$ as the linear restriction operator $Q$ fulfills $\mathbf{1}^TQ = \mathbf{1}$.
	
	For the first numerical tests, we choose $J=6, \lambda=[1.2,1.1,1.0,0.9,0.8,0.7]$, and $\mu=[1, 1, 1, 1, 1, 1]$ with three different capacities $k=8, 16$ and $32$. The dimension of each system is given by the formula $n=(k+1)^J$. On every grid we use an implementation of GMRES that exploits the tensor structure in matrix vector products and inner products as the smoother and perform $3$ pre- and postsmoothing steps for the first and second problem and $7$ pre- and postsmoothing steps for the third problem. We use the Moore--Penrose pseudo inverse to solve the singular coarse grid system. As initial guess we use the smallest singular vector from the matrix corresponding to the smallest grid and interpolate it up to the original dimension, which can be done in the setup at low cost in a bootstrap fashion~\cite{BrandtBrannickKahlLivshits2011}. In experiments we could observe that this typically leads to a better starting point than a random initial guess. The availability of this starting guess is an additional advantage of a multigrid approach in this context. Note that after every $V$-cycle we normalize our iterate to guarantee that $\mathbf{1}^Tx=1$. We start with a maximal rank of $30$ and increase the maximal rank by the factor $\sqrt{2}$ if after a $V$-cycle
	\begin{itemize}
		\item the iterate $X_i$ has reached its maximal rank and
		\item the residual norm of $X_i$ is not reduced or increased by $15\%$ percent compared to the residual norm of the iterate $X_{i-1}$.
	\end{itemize}

	We are able to increase the rank adaptively, which is also an advantage of our method, because the rank of the state vector is not known a priori. The table \ref{tab:overflow1} shows the number of $V$-cycles, the effective and maximal rank and the time needed to reduce the residual norm $\|Ax\|$ below $10^{-7}$ for these three different problems. 
	The effective rank of a TT-tensor $\mathcal{X}$ is the number $r_{\textit{eff}}$ such that $X$ requires the same amount of storage as a tensor with all TT-ranks equal to $r_{\textit{eff}}$, see~\cite{Savostyanov2012}.  
	
\begin{figure}[t]
		\centering
		\subfloat [][]{
			\includegraphics[width=.48\textwidth]{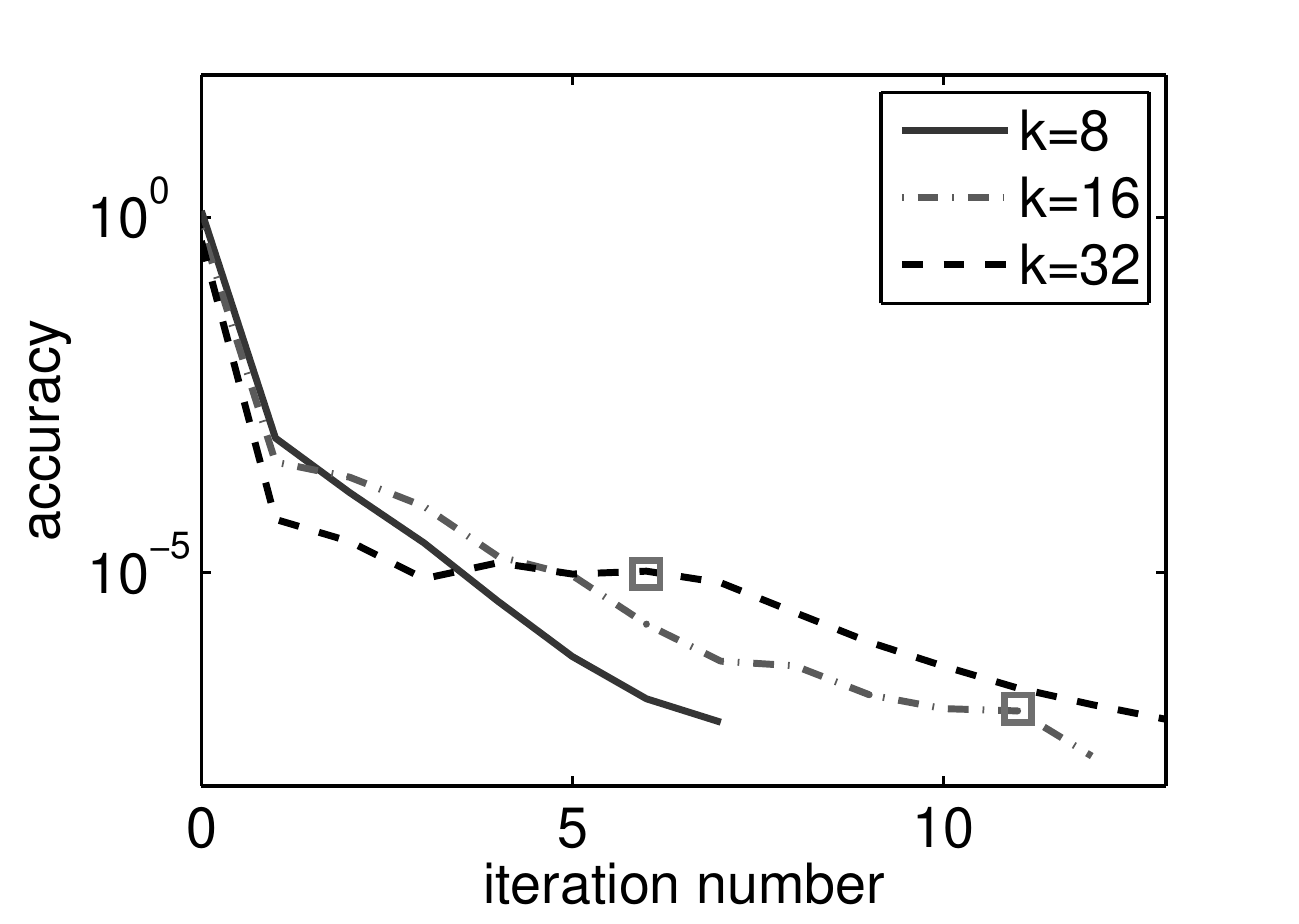}
			\label{fig:chJ6}
		}
		\hfill
		\subfloat[][]{
			\includegraphics[width=.48\textwidth]{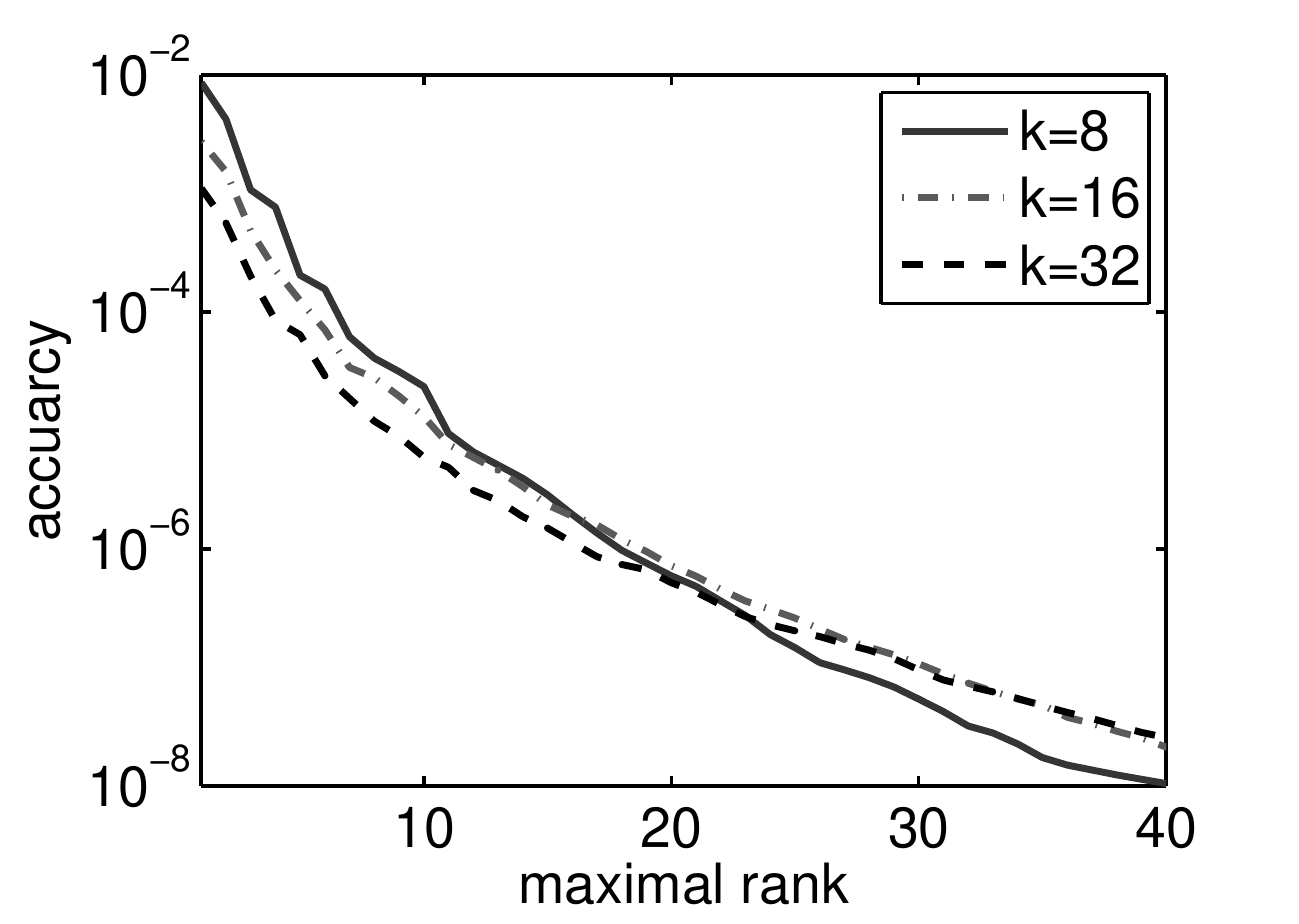}
			\label{fig:brJ6}
		}
		\label{fig:J6}
	\end{figure}

	The numerical tests show that for problems with $J=6$ the iteration number scales very nicely with the problem size.
	To guarantee that the used maximal rank in our tests is needed to approximate the solution and not only for the performance of our method, figure~\ref{fig:brJ6} shows the result of $\|AX\|$ where X is truncated which different truncation ranks. We see, that our method yields a solution with nearly "optimal" maximal rank. Figure~\ref{fig:chJ6} shows the convergence history of the three problems and the rectangles indicate in which iterations our method increases the maximal rank. 
	
	In order to compare our method with existing solution techniques we tested an implementation of alternating least squares with Tensor Train (ALS-TT)~\cite{KressnerMacedo2014} for problems with queue sizes $J=6$, $J = 7$, and $J = 8$. As ALS-TT requires an a priori rank-input we use the ranks determined by our multigrid method. The reduced problems which occur in ALS-TT were solved by the MATLAB backslash-solver. As a stopping criterion we chose the same accuracy requirement as for our multigrid method and aborted the iteration if the accuracy was not achieved after 4 hours. As can be seen in table~\ref{tab:overflow2} the absolute runtimes of ALS-TT are much worse than those of the multigrid algorithm, but due to the fact that this could be attributed to a non-optimal implementation of ALS-TT, we would rather like to point the attention to the scaling behaviour with respect to the queue size $k$ of both methods. Doubling the queue size and thus increasing the system size by a factor of $64$, the timings of the multigrid method grows by a factor of about $4$. ALS-TT on the other hand requires more than $100$ times as long. The slowing down of ALS-TT and other comparable methods like AMEn~\cite{KressnerMacedo2014} with growing model size has been attributed in the literature to the ill-conditioning of the occuring linear systems that need to be solved in these methods. The fact that the multigrid method reduces the ``local'' systems to small sizes and typically also yields better conditioned coarse representations of these operators might explain the drastically improved scaling behaviour. 
	The scaling of the two approaches with respect to the number of queues on the other hand is similar. Due to the fact that our coarsest system has the size $(n_i/2^{\#\text{levels}-1})^J$, i.e., it grows exponentially with $J$, we expect that using the pseudo-inverse to solve this system becomes infeasible for high dimensions $J$.
	
	For the remaining test problems we do not report comparisons with ALS-TT as we expect it to behave similarly.  However, note that  Table~\ref{tab:overflow2} shows some other test results for higher dimensions, for which the behavior is analogous to the problems with $J=6$.
	
	\begin{table}
		\begin{center}
			\begin{tabular}{|c|c|r|c|c|c|c|c|}\hline
				$J$&$k$&$n$&levels&iter&max.\ rank&eff.\ rank &time \\ \hline
				6 & 8 & $ 531,\!441$ &4& 7 & 30 & 19.8 & 4.8 \\
				6 & 16 & $24,\!137,\!569$ &5& 12 & 42 & 28.7 & 18.2 \\
				6 & 32 & $1,\!291,\!467,\!969$&6& 13 & 42 & 32.0 & 188.8\\ \hline
			\end{tabular}
		\end{center}
		\caption{Results of our multigrid algorithm for overflow queuing with 6 queues and varying capacity}
		\label{tab:overflow1}
	\end{table}
	
	\begin{table}
		\begin{center}
			\begin{tabular}{|c|c|r|c|c|}\hline
				$J$&$k$&$n$& multigrid time & ALS-TT time \\ \hline
				6 & 8 & $ 531,\!441$ &4.8 &    87.9 \\
				6 & 16 & $24,\!137,\!569$ &18.2&  $>$ 4 hours \\ \hline
				7 & 8 & $4,\!782,\!969$ & 7.8&  185.6 \\ \hline
				8 & 8 & $43,\!046,\!721$ & 11.9& 322.5  \\ \hline
			\end{tabular}
		\end{center}
		\caption{Comparison of timings between multigrid algorithm and ALS-TT}
		\label{tab:overflowvsals}
	\end{table}	
	
	\begin{table}
		\begin{center}
			\begin{tabular}{|c|c|r|c|c|c|c|c|}\hline
				$J$&$k$&$n$&levels&iter&max.\ rank&eff.\ rank &time \\ \hline
				7 & 8 & $4,\!782,\!969$ &4& 9 & 30 & 20.7 & 7.8 \\
				7 & 16 & $410,\!338,\!673$ &5& 15 & 30 & 21.2 & 27.6 \\ \hline
				8 & 8 & $43,\!046,\!721$ & 4 &10& 30 & 21.8 & 11.9 \\
				8 & 16 & $6,\!975,\!757,\!441$ &5& 17 & 30 & 22.1 & 40.1 \\ \hline
				9 & 8 & $387,\!420,\!489$ & 4&10& 30 & 22.76 & 16.0 \\
				9 & 16 & $118,\!587,\!876,\!497$ &5& 20 & 40 & 26.1 & 66.7 \\ \hline
				
			\end{tabular}
		\end{center}
		\caption{Results of our multigrid algorithm for overflow queuing with $\lambda=[1.2,1.1,1.0,0.9,0.8,0.7,\dots]$, and $\mu=[1, 1, 1, 1, 1, 1,\dots]$}
		\label{tab:overflow2}
	\end{table}
	
	\subsection{Manufacturing system with Kanban control}\label{subsec:kanban}
	The second example we consider is the so called Kanban model as it is discussed in~\cite{Buchholz1999}. 
	A finite number $J$ of cells is given and each of them contains a machine, a bulletin board and an output hopper. Parts have to go through these machines in a certain order. We assume that this order is given by the numeration of the machines. Only a certain number $k_i,\, i \dots, J$ of parts can enter a machine, which is controlled by the so called Kanban tickets, i.e., $k_i$ is the number of existing Kanban tickets for machine $i$. Each part which enters a machine gets a Kanban ticket and gives it back when it enters the subsequent machine. If no tickets are available in machine $i$ the part has to wait in the output hopper of the previous cell $i-1$ before entering the machine. The processing time and the time to move from one cell to the next are exponentially distributed with rate $\mu_i$ and rate $\omega_i$, respectively. As in the previous example we want to distinguish the local events and the synchronized ones. Note that the synchronized events exist only between neighboring machines. For sake of simplicity we initially do not consider the first and the last machine. 
	Each machine can be described by a Markov chain, where each state is characterized by three quantities:
	\begin{itemize}
		\item number of available tickets,
		\item number of parts being processed, 
		\item number of parts waiting for the next machine.
	\end{itemize}
	
	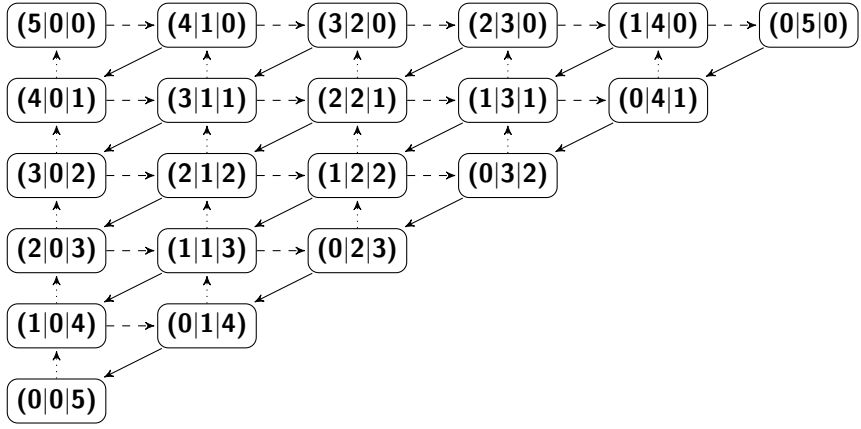
\begin{figure}
		\centering
		\begin{tikzpicture}[->,>=stealth',shorten >=1pt,auto,node distance=2cm ,main node/.style={rectangle, rounded corners, draw,font=\sffamily\bfseries}]
			
			\node[main node] (500) {(5$\mid$0$\mid$0)};
			\node[main node] (410) [right of=500] {(4$\mid$1$\mid$0)};
			\node[main node] (320) [right of=410] {(3$\mid$2$\mid$0)};
			\node[main node] (230) [right of=320] {(2$\mid$3$\mid$0)};
			\node[main node] (140) [right of=230] {(1$\mid$4$\mid$0)};
			\node[main node] (050) [right of=140] {(0$\mid$5$\mid$0)};
			
			\node[main node] (401) [below of=500,yshift=1cm] {(4$\mid$0$\mid$1)};
			\node[main node] (311) [below of=410,yshift=1cm] {(3$\mid$1$\mid$1)};
			\node[main node] (221) [below of=320,yshift=1cm] {(2$\mid$2$\mid$1)};
			\node[main node] (131) [below of=230,yshift=1cm] {(1$\mid$3$\mid$1)};
			\node[main node] (041) [below of=140,yshift=1cm] {(0$\mid$4$\mid$1)};
			
			\node[main node] (302) [below of=401,yshift=1cm] {(3$\mid$0$\mid$2)};
			\node[main node] (212) [below of=311,yshift=1cm] {(2$\mid$1$\mid$2)};
			\node[main node] (122) [below of=221,yshift=1cm] {(1$\mid$2$\mid$2)};
			\node[main node] (032) [below of=131,yshift=1cm] {(0$\mid$3$\mid$2)};
			
			\node[main node] (203) [below of=302,yshift=1cm] {(2$\mid$0$\mid$3)};
			\node[main node] (113) [below of=212,yshift=1cm] {(1$\mid$1$\mid$3)};
			\node[main node] (023) [below of=122,yshift=1cm] {(0$\mid$2$\mid$3)};
			
			\node[main node] (104) [below of=203,yshift=1cm] {(1$\mid$0$\mid$4)};
			\node[main node] (014) [below of=113,yshift=1cm] {(0$\mid$1$\mid$4)};
			
			\node[main node] (005) [below of=104,yshift=1cm] {(0$\mid$0$\mid$5)};
			
			\path[dotted,every node/.style={font=\sffamily\small}]
			(005) edge node { } (104)
			(104) edge node { } (203)  	
			(203) edge node { } (302)
			(302) edge node { } (401)
			(401) edge node { } (500)
			
			(014) edge node { } (113)  	
			(113) edge node { } (212)
			(212) edge node { } (311)
			(311) edge node { } (410)
			
			(023) edge node { } (122)
			(122) edge node { } (221)
			(221) edge node { } (320)
			
			(032) edge node { } (131)
			(131) edge node { } (230)
			
			(041) edge node { } (140);
			
			\path[dashed,every node/.style={font=\sffamily\small}]
			(500) edge node { } (410)
			(410) edge node { } (320)  	
			(320) edge node { } (230)
			(230) edge node { } (140)
			(140) edge node { } (050)
			
			(401) edge node { } (311)  	
			(311) edge node { } (221)
			(221) edge node { } (131)
			(131) edge node { } (041)
			
			(302) edge node { } (212)
			(212) edge node { } (122)
			(122) edge node { } (032)
			
			(203) edge node { } (113)
			(113) edge node { } (023)
			
			(104) edge node { } (014);	
			
			\path[every node/.style={font=\sffamily\small}]
			(410) edge node { } (401)
			(320) edge node { } (311)  	
			(230) edge node { } (221)
			(140) edge node { } (131)
			(050) edge node { } (041)
			
			(311) edge node { } (302)  	
			(221) edge node { } (212)
			(131) edge node { } (122)
			(041) edge node { } (032)
			
			(212) edge node { } (203)
			(122) edge node { } (113)
			(032) edge node { } (023)
			
			(113) edge node { } (104)
			(023) edge node { } (014)
			
			(014) edge node { } (005);
		\end{tikzpicture}
		\caption{Possible transitions between states of a machine with five Kanban tickets. Local transitions are depicted by solid arrows, synchronized transitions depending on the previous or subsequent machine by dashed or dotted arrows, respectively.}\label{fig:kanban_states}
	\end{figure}
	
	In Figure~\ref{fig:kanban_states} we illustrate the states of one machine with five tickets and characterize their transitions from one state to another as local or synchronized, where local transitions are the ones that are independent of the neighboring machines, i.e., the transition from the machine to the output hopper. 
	
	\begin{figure}
		\centering
		\begin{tikzpicture}[->,>=stealth',shorten >=1pt,auto,node distance=1cm ,main node/.style={rectangle, rounded corners, draw,font=\sffamily\bfseries}]
			
			\node[main node] (500a) {1};
			\node[main node] (410a) [right of=500a] {2};
			\node[main node] (320a) [right of=410a] {4};
			
			\node[main node] (401a) [below of=500a] {3};
			\node[main node] (311a) [below of=410a] {5};
			
			\node[main node] (302a) [below of=401a] {6};

			\node[main node] (500b) [right of=320a,xshift=1cm] {1};
			\node[main node] (410b) [right of=500b] {2};
			\node[main node] (320b) [right of=410b] {4};
			
			\node[main node] (401b) [below of=500b] {3};
			\node[main node] (311b) [below of=410b] {5};
			
			\node[main node] (302b) [below of=401b] {6};

			\node[main node] (500c) [right of=320b,xshift=1cm] {1};
			\node[main node] (410c) [right of=500c] {2};
			\node[main node] (320c) [right of=410c] {4};
			
			\node[main node] (401c) [below of=500c] {3};
			\node[main node] (311c) [below of=410c] {5};
			
			\node[main node] (302c) [below of=401c] {6};

			\path[dotted,every node/.style={font=\sffamily\small}]
			
			(302c) edge node { } (401c)
			(401c) edge node { } (500c)
			(311c) edge node { } (410c);

			\path[dashed,every node/.style={font=\sffamily\small}]
			(500b) edge node { } (410b)
			(410b) edge node { } (320b)  	
			(401b) edge node { } (311b); 	
			
			\path[every node/.style={font=\sffamily\small}]
			(410a) edge node { } (401a)
			(320a) edge node { } (311a) 	
			(311a) edge node { } (302a);
			
		\end{tikzpicture}
		
		\caption{Graphs corresponding to different transition types (left: local, middle: synchronized with previous machine, right; synchronized with subsequent machine) for a machine with two Kanban tickets. The states are ordered (and numbered) lexicographically, i.e., node 1 corresponds to $(2,0,0)$, node 2 corresponds to $(1,1,0)$ etc.}\label{fig:threegraphs}
	\end{figure}
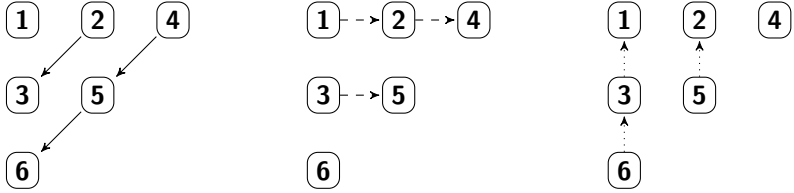
	
	Ordering the states lexicographically and distinguishing the three different types of transitions, we find three different graphs corresponding to one machine, as depicted in Figure~\ref{fig:threegraphs}. The matrices describing the different transition types are exactly the transpose of the adjacency matrices of these graphs, with all nonzero entries equal to $\mu_i$ (for local transitions), $\omega_{i-1}$ (for transitions depending on the previous machine) or $\omega_i$ (for transitions depending on the subsequent machine). In the local transition matrix, additionally the diagonal entries of nonzero columns are set to be $-\mu_i$. The corresponding matrices for a machine with two tickets are given as follows:
	\begin{equation}\label{eq:matrices2tickets}
		\begin{split}
			&E_i^{(i,i)} = \begin{pmatrix}
				0&0&0&0&0&0\\
				0&-\mu_i&0&0&0&0\\
				0&\mu_i&0&0&0&0\\
				0&0&0&-\mu_i&0&0\\
				0&0&0&\mu_i&-\mu_i&0\\
				0&0&0&0&\mu_i&0\\
			\end{pmatrix},\hspace{.05cm}
			E_i^{(i,i-1)}=\begin{pmatrix}
				0&0&0&0&0&0\\
				\omega_{i-1}&0&0&0&0&0\\
				0&0&0&0&0&0\\
				0&\omega_{i-1}&0&0&0&0\\
				0&0&\omega_{i-1}&0&0&0\\
				0&0&0&0&0&0\\
			\end{pmatrix},\\
			&\qquad\qquad\qquad\qquad
			E_i^{(i,i+1)} = \begin{pmatrix}
				0&0&\omega_i&0&0&0\\
				0&0&0&0&\omega_i&0\\
				0&0&0&0&0&\omega_i&\\
				0&0&0&0&0&0\\
				0&0&0&0&0&0\\
				0&0&0&0&0&0\\
			\end{pmatrix},\qquad i=2,\dots,J-1.
		\end{split}
	\end{equation}
	
	The difference of the first and the last machine to the others is that they have only one neighbor, so the number of states of the corresponding Markov chain is smaller. To be more precise, the last machine lacks an output hopper (i.e., as soon as a part is finished its ticket is available), whereas the first machine has no bulletin boards (i.e., in case there are tickets available, a new part immediately enters the machine). Thus both the first and the last machine can be described by only two of the three state quantities. 
	As such the states for the first and the last machine can also be read from figure~\ref{fig:kanban_states}. The first row of the state triangle describes the last machine and the first machine is given by the right-most diagonal.
	
	\begin{figure}[t]
		\centering
		\includegraphics[width=.6\textwidth]{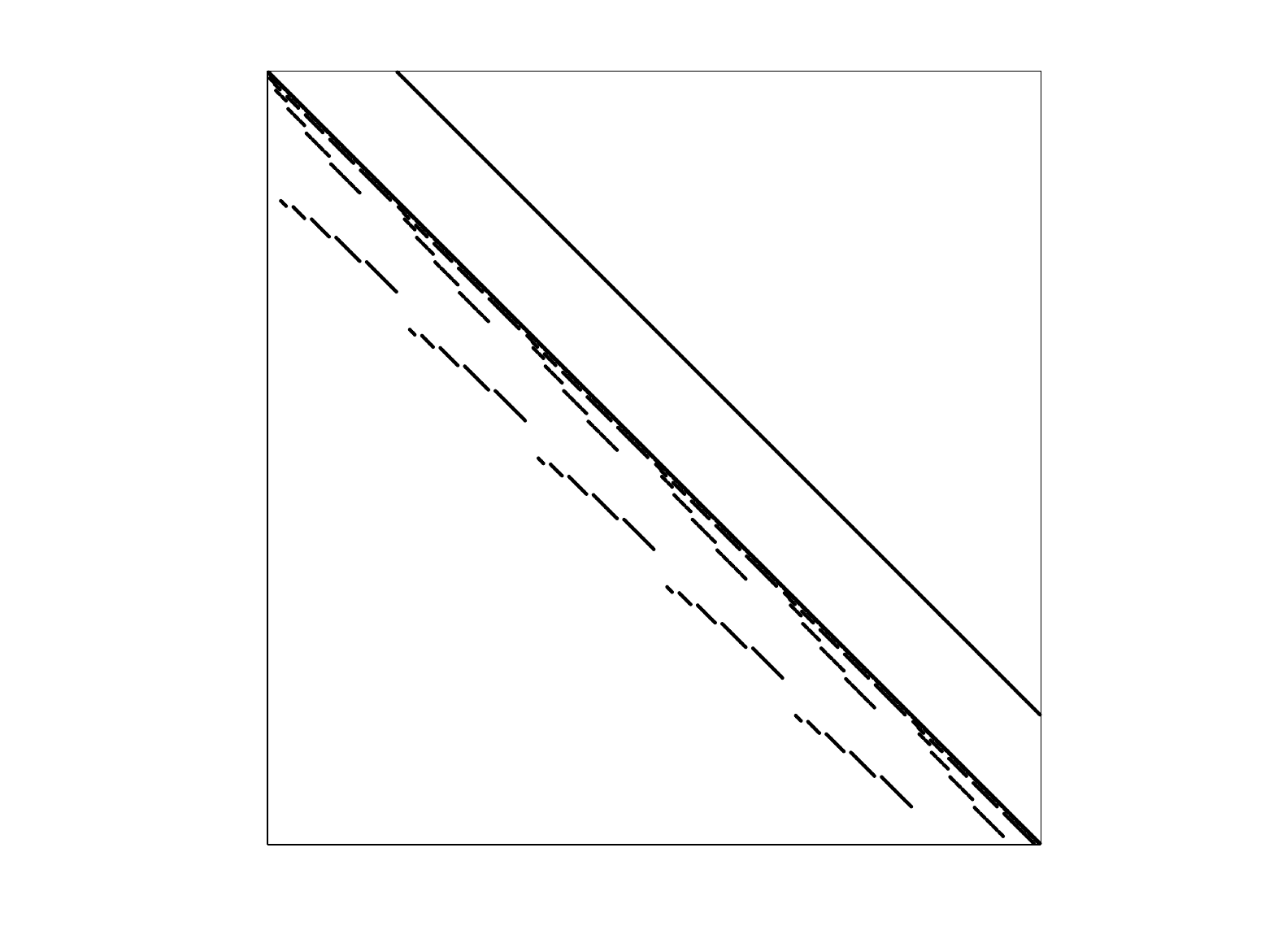}
		\caption{Nonzero structure of Kanban model for $J = 3, k = [5, 5, 5]$.}
		\label{fig:spykQ}
	\end{figure}
	
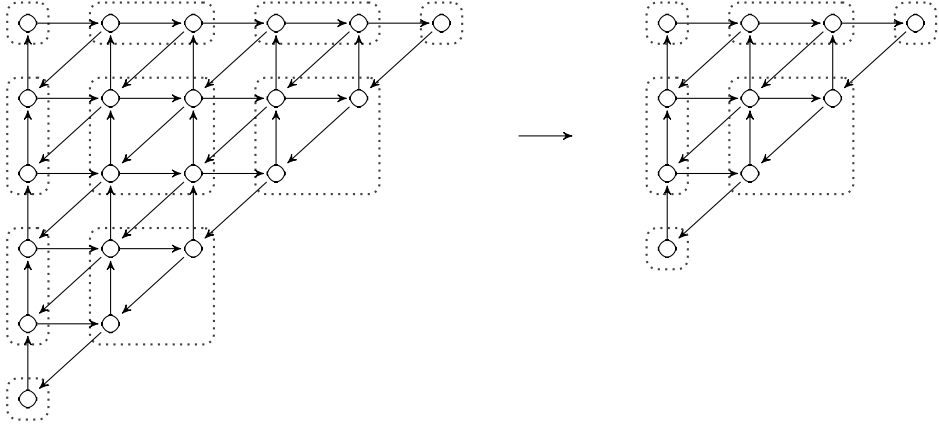
\begin{figure}
		\centering
		\begin{tikzpicture}[transform shape,->,>=stealth',shorten >=1pt,auto,node distance=1cm ,main node/.style={rectangle, rounded corners, draw,font=\sffamily\bfseries}]
			
			\node[main node] (500) {};
			\node[main node] (410) [right of=500, xshift=.1cm] {};
			\node[main node] (320) [right of=410, xshift=.1cm] {};
			\node[main node] (230) [right of=320,xshift=.1cm] {};
			\node[main node] (140) [right of=230,xshift=.1cm] {};
			\node[main node] (050) [right of=140,xshift=.1cm] {};
			
			\node[main node] (401) [below of=500] {};
			\node[main node] (311) [below of=410] {};
			\node[main node] (221) [below of=320] {};
			\node[main node] (131) [below of=230] {};
			\node[main node] (041) [below of=140] {};
			
			\node[main node] (302) [below of=401] {};
			\node[main node] (212) [below of=311] {};
			\node[main node] (122) [below of=221] {};
			\node[main node] (032) [below of=131] {};
			
			\node[main node] (203) [below of=302] {};
			\node[main node] (113) [below of=212] {};
			\node[main node] (023) [below of=122] {};
			
			\node[main node] (104) [below of=203] {};
			\node[main node] (014) [below of=113] {};
			
			\node[main node] (005) [below of=104] {};
			
			\node [rectangle, rounded corners, thick, draw=black!75!white, dotted, fit = (500), inner sep=0.15cm] {};
			\node [rectangle, rounded corners, thick, draw=black!75!white, dotted, fit = (050), inner sep=0.15cm] {};
			\node [rectangle, rounded corners, thick, draw=black!75!white, dotted, fit = (005), inner sep=0.15cm] {};
			
			\node [rectangle, rounded corners, thick, draw=black!75!white, dotted, fit = (401) (302), inner sep=0.15cm] {};
			\node [rectangle, rounded corners, thick, draw=black!75!white, dotted, fit = (203) (104), inner sep=0.15cm] {};
			\node [rectangle, rounded corners, thick, draw=black!75!white, dotted, fit = (410) (320), inner sep=0.15cm] {};
			\node [rectangle, rounded corners, thick, draw=black!75!white, dotted, fit = (230) (140), inner sep=0.15cm] {};
			
			\node [rectangle, rounded corners, thick, draw=black!75!white, dotted, fit = (311) (221) (212) (122), inner sep=0.15cm] {};
			
			\node [rectangle, rounded corners, thick, draw=black!75!white, dotted, fit = (131) (041) (032), inner sep=0.15cm] {};
			\node [rectangle, rounded corners, thick, draw=black!75!white, dotted, fit = (113) (023) (014), inner sep=0.15cm] {};
			
			\path[every node/.style={font=\sffamily\small}]
			(005) edge node { } (104)
			(104) edge node { } (203)  	
			(203) edge node { } (302)
			(302) edge node { } (401)
			(401) edge node { } (500)
			
			(014) edge node { } (113)  	
			(113) edge node { } (212)
			(212) edge node { } (311)
			(311) edge node { } (410)
			
			(023) edge node { } (122)
			(122) edge node { } (221)
			(221) edge node { } (320)
			
			(032) edge node { } (131)
			(131) edge node { } (230)
			
			(041) edge node { } (140);
			
			\path[every node/.style={font=\sffamily\small}]
			(500) edge node { } (410)
			(410) edge node { } (320)  	
			(320) edge node { } (230)
			(230) edge node { } (140)
			(140) edge node { } (050)
			
			(401) edge node { } (311)  	
			(311) edge node { } (221)
			(221) edge node { } (131)
			(131) edge node { } (041)
			
			(302) edge node { } (212)
			(212) edge node { } (122)
			(122) edge node { } (032)
			
			(203) edge node { } (113)
			(113) edge node { } (023)
			
			(104) edge node { } (014);	
			
			\path[every node/.style={font=\sffamily\small}]
			(410) edge node { } (401)
			(320) edge node { } (311)  	
			(230) edge node { } (221)
			(140) edge node { } (131)
			(050) edge node { } (041)
			
			(311) edge node { } (302)  	
			(221) edge node { } (212)
			(131) edge node { } (122)
			(041) edge node { } (032)
			
			(212) edge node { } (203)
			(122) edge node { } (113)
			(032) edge node { } (023)
			
			(113) edge node { } (104)
			(023) edge node { } (014)
			
			(014) edge node { } (005);

			\node[main node] (300c) [right of= 050, xshift = 2cm]{};
			\node[main node] (210c) [right of=300c, xshift=.1cm] {};
			\node[main node] (120c) [right of=210c, xshift=.1cm] {};
			\node[main node] (030c) [right of=120c,xshift=.1cm] {};
			
			\node[main node] (201c) [below of=300c] {};
			\node[main node] (111c) [below of=210c] {};
			\node[main node] (021c) [below of=120c] {};
			
			\node[main node] (102c) [below of=201c] {};
			\node[main node] (012c) [below of=111c] {};
			
			\node[main node] (003c) [below of=102c] {};
			
			\path[every node/.style={font=\sffamily\small}]
			(300c) edge node { } (210c)
			(210c) edge node { } (120c)  	
			(120c) edge node { } (030c)
			(201c) edge node { } (111c)  	
			(111c) edge node { } (021c)
			(102c) edge node { } (012c)
			
			(003c) edge node { } (102c)
			(102c) edge node { } (201c)
			(201c) edge node { } (300c)
			(012c) edge node { } (111c)
			(111c) edge node { } (210c)
			(021c) edge node { } (120c)
			
			(210c) edge node { } (201c)
			(120c) edge node { } (111c)
			(030c) edge node { } (021c)
			(111c) edge node { } (102c)
			(021c) edge node { } (012c)
			(012c) edge node { } (003c);
			
			\node [rectangle, rounded corners, thick, draw=black!75!white, dotted, fit = (300c), inner sep=0.15cm] {};
			\node [rectangle, rounded corners, thick, draw=black!75!white, dotted, fit = (030c), inner sep=0.15cm] {};
			\node [rectangle, rounded corners, thick, draw=black!75!white, dotted, fit = (003c), inner sep=0.15cm] {};
			
			\node [rectangle, rounded corners, thick, draw=black!75!white, dotted, fit = (210c) (120c), inner sep=0.15cm] {};
			\node [rectangle, rounded corners, thick, draw=black!75!white, dotted, fit = (201c) (102c), inner sep=0.15cm] {};
			\node [rectangle, rounded corners, thick, draw=black!75!white, dotted, fit = (111c) (021c) (012c), inner sep=0.15cm] {};	
			
			\node (1) [right of=041, xshift=1cm,yshift=-0.5cm] {};
			\node (2) [right of=1] {};
			\path[every node/.style={font=\sffamily\small}]
			(1) edge node { } (2);
		\end{tikzpicture}
		\caption{Choice of aggregates in a three-level method for a machine with five tickets.}\label{fig:kanban_aggregates}
	\end{figure}

	The sparsity structure of the whole system with three machines and five tickets each is given in figure~\ref{fig:spykQ}. In contrast to the overflow queuing model the graphs corresponding to the local matrices have some unreachable states, thus full coarsening with these matrices is not possible and in addition the geometric structure of these matrices does not bear a likeness to the structure of the system matrix depicted in figure~\ref{fig:spykQ}, for example they have some zero columns. We thus base our coarsening approach on the accumulated local structure as shown in figure~\ref{fig:kanban_states} and use a simple aggregation approach with constant basis function. On the left of figure~\ref{fig:kanban_aggregates} the selected aggregates for a problem with five tickets are shown. For the first and last machine the aggregates are chosen accordingly based on the identification with the right-most diagonal and first row, respectively. 
	By enumerating the aggregates along the diagonals of the grid (starting in the upper left corner, just like the original numbering of the nodes in Figure~\ref{fig:threegraphs}) we preserve the connection structure of figure~\ref{fig:kanban_states} and thus are able to continue aggregating according to this geometry to end up with a multilevel approach. On the right of figure~\ref{fig:kanban_aggregates} the choice of aggregates for further coarsening to a third level are shown. Note that the number of tickets in one machine should be $2^k+1$ to allow for such an aggregation and we limit our analysis to such Kanban systems. This limitation is simply a technical issue and we suspect that irregular numbers of tickets can be dealt with by adapting the aggregation accordingly.
	
	If the number of tickets is the same on every machine, which is denoted by $k$, the size of the system matrix is given by the formula $n=(k+1)^2 \cdot ((k+1)(k+2)/2)^{(J-2)}$.
	With the aggregation as described before we are able to reduce our problem to a size of 
	$3^2\cdot6^{(J-2)}$. To reduced this to $2^2\cdot3^{(J-2)}$, we aggregate state $1$ and $2$,
	state $4$ and $5$ and state $3$ and $6$ (numbering as in figure~\ref{fig:threegraphs}) on the second to last level.
	
	For our numerical tests we choose the aggregation approach as described before, therefore the choice of interpolation and restriction operator is given by standard aggregation interpolation (the restriction operator is just the transpose of the interpolation operator), see, e.g.,~\cite{Vanek1995}, and the coarse grid matrix $A_c$ is the corresponding (Petrov--)Galerkin operator as before. Note that the column sum of $A_c$ is equal zero, because every column in the restriction matrix has only one entry with the value one.
	As a smoother, we use one approximate Gau{\ss}-Seidel iteration on each level, see~\cite{Greenbaum1997,TrottenbergOsterleeSchueller2001}. The lower and strictly triangular part $L$ and $U$ of our system matrix $A$ is given in a Tensor representation, similar to A, see~\cite{UysalDayar1998}. To apply the action of $L^{-1}$, we use the AMEn iteration from ~\cite{DolgovSavostyanov2013a,DolgovSavostyanov2013b} (with stopping accuracy $10^{-7}$ and a maximum of 20 sweeps) as implemented in the TT-Toolbox. Fortunately this turns out to be much easier than applying AMEn to the original problem. The reason may be that the matrix $L$ is non-singular in contrast to $A$. Note that the linear systems that occur in AMEn for $L$ are much better conditioned than the same systems for $A$.
	In addition, the structure of the Kanban control model, which is similar to a flow problem, also suggests that Gau{\ss}-Seidel is a good choice as a smoother, as it is known to be efficient for those problems.
	
	The other ingredients of the method are the same as in example~\ref{subsec:overflow}. Table~\ref{tab:kanban} shows the results of some numerical tests for different numbers of machines and tickets. We can observe that the method works well for all tested problems. Again we see that a higher rank is needed with increasing dimensions of the problem. Apart from the rank increase we can again observe that the method scales reasonably well with growing problem size. The time the method takes is much higher than for the overflow test problems (mainly due to the use of the Gau{\ss}-Seidel smoother and the larger coarse grid systems) but this has to be expected as the Kanban problem is also much harder to solve for other solvers, see, e.g.,~\cite{KressnerMacedo2014}. 
	
	\begin{table}
		\begin{center}
			\begin{tabular}{|c|c|r|c|c|c|c|c|}\hline
				$J$&$k$&$n$&levels&iter&max.\ rank&eff.\ rank &time \\ \hline
				6 & 3 & $160,\!000$ &3& 11 & 22 & 11.4 & 308.8\\
				6 & 5 & $7,\!001,\!316$ &4& 14 & 41 & 17.5 & 1136.5\\ 
				6 & 9 & $915,\!062,\!500$ &5& 23 & 59 & 29.6 & 12526.1\\ \hline
				7 & 3 & $1,\!600,\!000$ &3& 14 & 30 & 15.2 & 689.8\\
				7 & 5 & $147,\!027,\!636$ &4& 16 & 59 & 28.0 & 2051.1\\ \hline
				8 & 3 & $16,\!000,\!000$ &3& 16 & 42 & 21.9 & 1264.4\\
				8 & 5 & $3,\!087,\!580,\!356$ &4& 15 & 59 & 33.8 & 2783.8\\ \hline	
			\end{tabular}
		\end{center}
		\caption{Results of our multigrid algorithm for Kanban model with varying number of machines and tickets and $\mu_i = 1, \omega_i = .1$ for all $i$.}
		\label{tab:kanban}
	\end{table}

	\section{Conclusion}\label{sec:conclusion}
	We discussed how to develop multigrid methods for tensor-structured problems, more specifically for structured Markov chains, which on the one hand exploit the tensor structure for efficient matrix and vector operations and at the same time keep this structure intact across all grids. We illustrated the behavior of our method by investigating two standard model problems.
	
	Topics for future research include exploring different coarsening strategies, e.g., aggregation across submodels (i.e., queues or machines for the overflow queuing or Kanban system, respectively) for large scale problems.
	
	A further extension of this work is generalizing the structure-preserving multigrid approach to other classes of structured Markov chains, e.g., multi-server multi-queue models.
	
	\bibliographystyle{siam}
	\bibliography{markov}{}
\end{document}